\documentclass[11pt,oneside]{article}
\usepackage{amssymb}
\usepackage{amsmath}
\usepackage{latexsym}
\usepackage{amsfonts}
\usepackage{amscd}
\usepackage{amsthm}
\usepackage{fancyhdr}

\usepackage[headings]{fullpage}

\newcommand{\gt}[1]{\mathfrak{#1}}

\newcommand{\mc}[1]{\mathcal{#1}}

\newcommand{\comment}[1]{}

\newcommand{\RR}{{\mathbb R}}
\newcommand{\CC}{{\mathbb C}}
\newcommand{\ZZ}{{\mathbb Z}}
\newcommand{\ZZh}{\ZZ+\tfrac{1}{2}}
\newcommand{\hZZ}{\tfrac{1}{2}\ZZ}

\newcommand{\hh}{{\mathbf h}}
\newcommand{\FF}{{\mathbb F}}

\newcommand{\fset}{{\Sigma}}
\newcommand{\hfset}{{\Delta}}    

\newcommand{\ii}{{\bf i}}
\newcommand{\lab}{{\langle}}    
\newcommand{\llab}{{\langle\!\langle}}  
\newcommand{\rab}{{\rangle}}    

\newcommand{\ldb}{(\!(}
\newcommand{\rdb}{)\!)}

\newcommand{\qj}{{\bf j}}
\newcommand{\qk}{{\bf k}}
\newcommand{\ql}{{\bf l}}

\newcommand{\LL}{{\Lambda}}

\newcommand{\End}{\operatorname{End}}
\newcommand{\Aut}{\operatorname{Aut}}
\newcommand{\Id}{\operatorname{Id}}
\newcommand{\tr}{\operatorname{{\sf tr}|}}

\newcommand{\Res}{\operatorname{Res}}
\newcommand{\Sing}{\operatorname{Sing}}
\newcommand{\Lie}{\operatorname{Lie}}
\newcommand{\Img}{\operatorname{Im}}
\newcommand{\Span}{\operatorname{Span}}
\newcommand{\rank}{\operatorname{rank}}

\newcommand{\GL}{\operatorname{\textsl{GL}}}      
\newcommand{\PSL}{\operatorname{\textsl{PSL}}}    
\newcommand{\SL}{\operatorname{\textsl{SL}}}      
\newcommand{\SO}{\operatorname{\textsl{SO}}}    
\newcommand{\Mp}{\operatorname{\textsl{Mp}}}    
\newcommand{\Sy}{\operatorname{\textsl{Sp}}}  
\newcommand{\U}{\operatorname{\textsl{U}}}        
\newcommand{\SU}{\operatorname{\textsl{SU}}}    

\newcommand{\MM}{\mathbb{M}}    

\newcommand{\Ru}{\operatorname{\textsl{Ru}}}

\newcommand{\Wy}{\operatorname{Weyl}}

\newcommand{\Wyx}{\widetilde{\operatorname{Weyl}}}

\newcommand{\WS}{\forall}            

\newcommand{\vn}{V^{\natural}}     

\newcommand{\aru}{{A_{{\Ru}}}}  
\newcommand{\arutd}{{\tilde{A}_{{\Ru}}}}

\newcommand{\wru}{{\tilde{\forall}_{{\Ru}}}} 

\newcommand{\wu}{\WS(\gt{u})}

\newcommand{\cgset}{\Omega}   
\newcommand{\cgU}{\cgset_{\U}}

\newcommand{\wcgset}{\mho}    
\newcommand{\wcgU}{\wcgset_{\U}}
\newcommand{\wcgRu}{\tilde{\wcgset}_{\Ru}}

\newcommand{\cfalg}{\mc{A}}   

\newcommand{\cej}{\jmath}

\newcommand{\ceJ}{J}

\newcommand{\ossp}{\delta}

\newcommand{\ogi}{\theta}   

\newcommand{\ogih}{\ogi^{1/2}}

\newcommand{\vac}{{\bf 1}}      
\newcommand{\cas}{{\bf \omega}} 

\newtheorem{thm}{Theorem}[section]
\newtheorem{cor}[thm]{Corollary}

\newtheorem{prop}[thm]{Proposition}

\theoremstyle{definition}
\newtheorem*{defn}{Definition}

\theoremstyle{remark}
\newtheorem*{rmk}{Remark}

\numberwithin{equation}{subsection}
\newtheorem*{eg}{Example}

\theoremstyle{plain}
\newtheorem*{sthm}{Theorem}

\begin{document}

\title{
    \textsc{Moonshine for Rudvalis's sporadic group II}
    \footnote{{\it Mathematics Subject Classification (2000)}
                17B69, 
                20D08.
                }
          }

\author{John F. Duncan
          \footnote{
          Harvard University,
          Department of Mathematics,
          One Oxford Street,
          Cambridge, MA 02138,
          U.S.A.}
          {}\footnote{
          Email: {\tt duncan@math.harvard.edu};\;
          homepage: {\tt
          http://www.math.harvard.edu/\~{}jfd/}
               }
               }

\date{October 18, 2008}

\setcounter{section}{-1}


\maketitle

\begin{abstract}
In Part I we introduced the notion of enhanced vertex operator
superalgebra, and constructed an example which is self-dual, has
rank $28$, and whose full symmetry group is a seven-fold cover of
the sporadic simple group of Rudvalis. In this article we construct
a second enhanced vertex operator superalgebra whose full
automorphism group is a cyclic cover of the Rudvalis group. This new
example is self-dual and has rank $-28$. As in Part I, we can
compute all the McKay--Thompson series associated to the action of
the Rudvalis group explicitly. We observe that these series, when
considered together with those of Part I, satisfy a genus zero
property.
\end{abstract}

\tableofcontents

\section{Introduction}\label{sec:intro}

Monstrous Moonshine is, in a word, the unexpected connection between
modular functions\footnote{Here, by {\em modular function} we
understand a holomorphic function on the upper half plane that is
invariant for the action of some discrete subgroup of ${\rm
PSL}_2(\mathbb{R})$ that is commensurable with the modular group
${\rm PSL}_2(\mathbb{Z})$.}, and the largest sporadic simple group,
the Fischer--Griess Monster. Since the initial observations of McKay
and Thompson
\begin{gather}\label{eqn:intro_MTObs}
     \begin{split}
     196884&=1+196883\\
     21493760&=1+196883+21296876\\
     864299970&=(2)1+(2)196883+21296876+842609326\\
     &\&c.
     \end{split}
\end{gather}
relating coefficients of the Fourier expansion of Klein's modular
invariant (on the left) with degrees of irreducible representations
of the Monster group (on the right), much fascinating and beautiful
mathematics has been developed in order to explicate these
coincidences.

It is very striking that the largest sporadic simple group has a
connection with modular functions. Even more surprising are the
particular special properties that these ``Monstrous'' modular
functions $F_g(\tau)$ for $g\in\mathbb{M}$ appear to satisfy.
Roughly speaking, each one is a generator for the field of functions
on the curve $({\mathbb{H}/\Gamma_g})^*$ defined by it's invariance
group $\Gamma_g\subset{\rm PSL}_2(\mathbb{R})$, and in particular,
the curve in each case has genus zero. Armed only with the character
table of the (at that time conjectural) Monster group, and with the
observations (\ref{eqn:intro_MTObs}) of McKay and Thompson, Conway
and Norton were able to collect a wealth of information about what
the connection should entail \cite{ConNorMM}, and in particular,
formulated their {\em Moonshine conjectures:}
\begin{quote}
There exists a graded vector space $V=\bigoplus_n V_n$ with an
action by $\MM$ such that for each $g\in \MM$ the graded trace
function $F_g(\tau)=\sum_n\tr_{V_n}gq^n$ is a generator for the
function field determined by a discrete group $\Gamma_g<\PSL_2(\RR)$
of genus zero.
\end{quote}
The proof of these conjectures is due to Borcherds \cite{BorMM}, and
is a fine example of the far reaching developments in mathematics
that have occurred since the inception of Moonshine. Important roles
in the proof of the Moonshine conjectures are played by {\em vertex
algebras} introduced by Borcherds \cite{BorPNAS}, and by a
particular {\em vertex operator algebra (VOA)} constructed by
Frenkel--Lepowsky--Meurman \cite{FLM} called the {\em Moonshine VOA}
and denoted $\vn$. Indeed, the VOA $\vn$ furnishes an algebraic
structure whose full automorphism group is the Monster group $\MM$
and whose underlying vector space is a natural candidate for the $V$
in the statement of the Moonshine conjectures.

In brief, a vertex operator algebra is a vertex algebra equipped
with extra structure (and conditions). A simple theme in this
article and its companion \cite{DunVARuI} (see also \cite{Dun_VACo})
is to explore what happens when one specializes further the
structure of vertex algebra; we arrive at the notion of {\em
enhanced vertex algebra}.

Having observed aspects of Monstrous Moonshine, it is natural to ask
if there are analogous phenomena with other finite simple groups in
place of the Monster. Perhaps of greatest importance are the
sporadic simple groups, since at present they resist a uniform
description, and must be studied via somewhat ad. hoc. methods. Any
approach that might furnish a uniform treatment of the sporadic
groups would be of great interest. Curiously, many (indeed most) of
the $26$ sporadic simple groups are {\em involved} in the
Monster\footnote{A group is {\em involved} in the Monster if it can
be identified with a quotient of some subgroup of the Monster.}, and
for these sporadic groups there is substantial evidence of Moonshine
type phenomena (c.f. \cite{QueMFnsFGps}). In fact, these
observations for sporadic groups that are involved in the Monster
({\em monstrous sporadics}) may be regarded as the content of the
{\em generalized Moonshine conjectures} due to Norton:
\begin{quote}
For each $g\in\MM$ there is graded vector space $V_g=\bigoplus_n
(V_g)_n$ with a (projective) action by $C_{\MM}(g)$ such that for
each $h\in C_{\MM}(g)$ the graded trace function
$Z(g,h,\tau)=\sum_n{\sf tr}|_{(V_g)_n}hq^n$ is either a generator
for the function field determined by a discrete group
$\Gamma_{g,h}<{\PSL}_2(\RR)$ of genus zero, or is constant.
\end{quote}
In the preceding article \cite{DunVARuI} we focused attention on a
particular sporadic group that is not involved in the Monster: the
sporadic simple group of Rudvalis. We introduced the notion of {\em
enhanced vertex operator superalgebra} which is, in brief, a vertex
operator superalgebra equipped with extra structure. We constructed
a particular example $\aru$ whose full automorphism group is a
cyclic cover of the Rudvalis group. Even more than this, the object
$\aru$ is somewhat distinguished among vertex algebras in that it is
self-dual (that is, has a unique simple module) and satisfies a
strong vanishing condition on certain homogeneous subspaces of low
degree. This suggests that $\aru$ may admit a characterisation
analogous to those which hold for classical objects such as a the
Golay code, and the Leech lattice (c.f. \cite{Con69}), and which
have been established to some extent for other distinguished
(enhanced) vertex operator algebras such as the Moonshine VOA $\vn$
(c.f. \cite{FLM}, \cite{DonGriLam_UniqMnVOA},
\cite{LamYam_ChznMoonshineVOA}, and see also \cite{Dun_VACo}). The
action of the Rudvalis group $\Ru$ on $\aru$ allows one to attach
modular functions (and modular forms) to elements $g\in \Ru$, and in
fact, the extra structure that constitutes the particular enhanced
vertex operator superalgebra structure on $\aru$ gives rise to
Jacobi forms in a natural way, and we obtain, for example, the
following analogues of McKay and Thompson's observations
(\ref{eqn:intro_MTObs})
\begin{gather}\label{eqn:intro_RuIObs}
     \begin{split}
     378&=378\\
     784&=1+783\\
     20475&=20475\\
     92512&=(2)378+406+91350\\
          &\&c.
     \end{split}
\end{gather}
relating coefficients of a particular Jacobi form (on the left) with
degrees of irreducible representations of the Rudvalis group (on the
right).

In the present article we construct a second enhanced vertex
operator algebra $\wru$ whose full automorphism group is a cyclic
cover of the sporadic group of Rudvalis. The object $\wru$ is in
some sense a mirror to $\aru$; it shares many properties with $\aru$
in a kind of ``reciprocal'' way. For example, $\wru$ is self-dual,
as is $\aru$, and $\wru$ has rank $-28$ --- just the opposite of the
rank of $\aru$. Even at a practical level, $\wru$ admits a
construction which runs directly parallel to that given for $\aru$
in the preceding article \cite{DunVARuI}, and the format of the
present article has been organized so as to emphasize this.

The structure on $\wru$ gives rise to Jacobi forms just as in the
case of $\aru$, and we obtain a second collection of coincidences in
analogy with (\ref{eqn:intro_MTObs}), (\ref{eqn:intro_RuIObs}).
\begin{gather}\label{eqn:intro_RuIIObs}
     \begin{split}
     406&=406\\
     784&=1+783\\
     31465&=1+783+3276+27405\\
     114464&=(2)378+(3)406+3654+45500+63336\\
     &\&c.
     \end{split}
\end{gather}
This may reassure the reader who noticed an absence of some small
representations of $\Ru$ in (\ref{eqn:intro_RuIObs}); relations
(\ref{eqn:intro_RuIObs}) and (\ref{eqn:intro_RuIIObs}) incorporate
all of the irreducible characters of $\Ru$ of degree less than
$10^5$.

Each of the enhanced vertex operator superalgebras $\aru$ and $\wru$
admit canonical automorphisms of order two. (In the former case the
automorphism is that associated to the super structure --- a natural
$\ZZ/2$-grading on the underlying vector space. In the latter case
there is again a natural $\ZZ/2$-grading, but the chosen vertex
algebra structure is nonetheless purely even.) Since each vertex
algebra is self-dual, there are unique canonically twisted modules
$(\aru)_{\ogi}$ and $(\wru)_{\ogi}$ over $\aru$ and $\wru$,
respectively, and the group $\Ru$ acts projectively on each module.
This allows one to attach further functions with modular properties
to elements $g\in \Ru$, and in the last part of the article we use
these four spaces: $\aru$, $\wru$, and their canonically twisted
modules, to construct certain distinguished functions
$\tilde{F}^{A}_{g}(\tau)$ and $\tilde{F}^{\forall}_{g}(\tau)$ for
each $g\in \Ru$.

Our main motivation for introducing this second realization of the
Rudvalis group by vertex operators on $\wru$ is the following
observation.
\begin{quote}
For each $g\in \Ru$ the functions $\tilde{F}^A_g(\tau)$ and
$\tilde{F}^{\forall}_g(\tau)$ span a (two dimensional)
representation of a discrete subgroup of $\PSL_2(\RR)$ that is
commensurable with $\PSL_2(\ZZ)$ and has genus zero.
\end{quote}
Thus we obtain a genus zero property for the sporadic group of
Rudvalis.

As in the case of the Monster, the reason for this property remains
mysterious. We do have the advantage here that the actions of $\Ru$
are sufficiently transparent that it is not hard to compute the
functions $\tilde{F}^A_g$ and $\tilde{F}^{\forall}_g$ explicitly,
and check the invariance groups directly. The genus zero property
for $\Ru$ is certainly weaker than that which holds for the Monster
group. For example, the functions $\tilde{F}^A_g$ and
$\tilde{F}^{\forall}_g$ are constant for some choices $g\in \Ru$
(giving rise to the trivial two dimensional representation of
$\PSL_{2}(\ZZ)$), but this should be compared with the statement of
the generalized Moonshine conjectures above.

Evidently the genus zero phenomena, which is such a surprising
aspect of Monstrous Moonshine, extends in some sense to at least one
sporadic group beyond the Monster. It is now of great importance to
determine if this behavior is not shared by the other non-monstrous
sporadics. 

\medskip

The main result of this article is the following theorem.
\begin{sthm}[\ref{thm:Ru:main}]
The quadruple $(\wru,Y,\vac,\wcgRu)$ is a self-dual enhanced VOA of
rank $-28$. The full automorphism group of $(\wru,\wcgRu)$ is a
four-fold cover of the sporadic simple group of Rudvalis.
\end{sthm}

\subsection{Outline}\label{sec:intro:outline}

The format of this article runs parallel to that of \cite{DunVARuI},
with Weyl algebras playing the role formerly taken by Clifford
algebra. For this reason, we recall only briefly the background
material that is discussed more fully in \cite{DunVARuI}.

In \S\ref{sec:VOAs} we review some facts about VOAs\footnote{Here,
and from here on, we will usually suppress the ``super'' in
super-objects, so that unless extra clarification is necessary,
superspaces and superalgebras will be referred to as spaces and
algebras, respectively, and the term VOA for example, will be used
even when the underlying vector space comes equipped with a
$\ZZ/2$-grading.}. We recall the definition in
\S\ref{sec:VOAs:struc}, and we recall the notion of vertex Lie
algebra in \S\ref{sec:VOAs:VLAs}. The notion of enhanced VOA is
reviewed in \S\ref{sec:eVOAs}.

In \S\ref{sec:weylalgs:struc} we make our conventions regarding Weyl
algebras, and in \S\ref{sec:weylalgs:meta} the metaplectic
groups. 
In \S\ref{sec:weylalgs:VOAs} we review a previously known
construction of VOA structure on certain Weyl algebra modules
constructed in turn from a finite dimensional vector space with
antisymmetric bilinear form, and in \S\ref{sec:weylalgs:Herm} we
make some conventions regarding the Hermitian structures arising in
the case that the initial vector space comes equipped also with a
suitable Hermitian form.

In \S\ref{sec:GLn} we present a family of enhanced VOAs whose full
automorphism groups are the general linear groups $\GL_N(\CC)$. This
section prepares the way for our main example, which is considered
in \S\ref{sec:Ru}. The construction is parallel to that of the main
example of Part I \cite{DunVARuI}, although it is less involved,
since the enhanced conformal structure admits an extremely simple
description in terms of the Conway--Wales lattice. We describe this
construction in \S\ref{sec:Ru:geom}, and then the construction of
$\wru$ is given in \S\ref{sec:Ru:const}. We determine the symmetry
group of $\wru$ in \S\ref{sec:Ru:symms}.

In \S\ref{sec:series} we furnish explicit expressions for the
McKay--Thompson series associated to the Rudvalis group via its
action on $\wru$. The results of this section, along with those of
the corresponding section of Part I allow us to formulate a genus
zero property for certain modular functions associated to the
Rudvalis group. This property is described in
\S\ref{sec:series:MbeyM}. We also mention in
\S\ref{sec:series:MbeyM} an extension of this to the third group of
Janko.


\subsection{Notation}\label{sec:intro:notation}

Our notation follows that of \cite{DunVARuI}. For example, $\ii$
denotes a square root of $-1$ in $\CC$. We use $\FF^{\times}$ to
denote the non-zero elements of a field $\FF$. More generally,
$A^{\times}$ shall denote the set of invertible elements in an
algebra $A$. In this article we shall use $\FF$ to denote either
$\RR$ or $\CC$.

A {\em superspace} is a vector space with a grading by
$\ZZ/2=\left\{\bar{0},\bar{1}\right\}$. When $M$ is a superspace, we
write $M=M_{\bar{0}}\oplus M_{\bar{1}}$ for the superspace
decomposition, and for $u\in M$ we set
$|u|=\gamma\in\{\bar{0},\bar{1}\}$ when $u$ is $\ZZ/2$-homogeneous
and $u\in M_{\gamma}$. The dual space $M^*$ has a natural superspace
structure such that $(M^*)_{\gamma}=(M_{\gamma})^*$ for $\gamma\in
\{\bar{0},\bar{1}\}$. The space $\End(M)$ admits a structure of Lie
superalgebra when equipped with the {\em Lie superbracket}
$[\cdot\,,\cdot]$ which is defined so that
$[a,b]=ab-(-1)^{|a||b|}ba$ for $\ZZ/2$-homogeneous $a,b$ in
$\End(M)$. All formal variables will be regarded as even, so that
$M\ldb z\rdb_{\gamma}=M_{\gamma}\ldb z\rdb$ for example, for
$\gamma\in \ZZ/2$.

We write $\bigwedge(M)$ for the full exterior algebra of a vector
space $M$. We write
$\bigwedge(M)=\bigwedge(M)^0\oplus\bigwedge(M)^1$ for the parity
decomposition of $\bigwedge(M)$, and we write
$\bigwedge(M)=\bigoplus_{k\geq 0}\bigwedge^k(M)$ for the natural
$\ZZ$-grading on $\bigwedge(M)$.

We denote by $D_z$ the operator on formal power series which is
formal differentiation in the variable $z$, so that if $f(z)=\sum
f_nz^{-n-1}\in V[[z^{\pm 1}]]$ is a formal power series with
coefficients in some space $V$, we have
$D_zf(z)=\sum_n(-n)f_{n-1}z^{-n-1}$. For $m$ a non-negative integer,
we set $D_z^{(m)}=\tfrac{1}{m!}D_z^m$.

We use $\eta(\tau)$ to denote the Dedekind eta function.
\begin{gather}\label{eqn:intro:notation:eta}
        \eta(\tau)=q^{1/24}\prod_{n\geq 1}(1-q^n)
\end{gather}
Here $q=e^{2\pi\ii\tau}$ and $\tau$ is a variable in the upper half
plane, which we denote by $\hh$. Recall also the Jacobi theta
function $\vartheta(z|\tau)$ which we normalize so that
\begin{gather}\label{eqn:intro:notation:theta}
    \vartheta(z|\tau)=\sum_{m\in\ZZ}
        e^{2\pi\ii z m+\pi\ii\tau m^2}
\end{gather}
for $\tau\in\hh$ and $z\in\CC$. According to the Jacobi Triple
Product Identity we have
\begin{gather}\label{eqn:intro:tripprod}
    \vartheta(z|\tau)
    =\prod_{m\geq 0}(1-q^{m+1})
        (1+e^{2\pi\ii z}q^{m+1/2})
        (1+e^{-2\pi\ii z}q^{m+1/2})
\end{gather}
with $q=e^{2\pi\ii\tau}$ as before.

\medskip

The most specialized notations arise in \S\ref{sec:weylalgs}. We
include here a list of them, with the relevant subsections indicated
in brackets. They are grouped roughly according to similarity of
appearance, rather than by order of appearance, so that the list may
be easier to search through, whenever the need might arise.
\begin{small}
\begin{list}{}{     \itemsep -3pt
                    \topsep 3pt
                         }
\item[$\gt{a}$]     A complex vector space with non-degenerate
Hermitian form (\S\ref{sec:weylalgs:Herm}).

\item[$\gt{a}^*$]   The dual space to $\gt{a}$
(\S\ref{sec:weylalgs:Herm}).

\item[$\gt{u}$]  A real or complex vector space of even dimension with
non-degenerate bilinear form. In the case that $\gt{u}=\gt{a}\oplus
\gt{a}^*$, the bilinear form is assumed to be $1/2$ times the
symmetric linear extension of the natural pairing between $\gt{a}$
and $\gt{a}^*$ (\S\ref{sec:weylalgs:Herm}).

\item[$\{a_i\}_{i\in\hfset}$]  A basis for $\gt{a}$, orthonormal in
the sense that $(a_i,a_j)=\delta_{ij}$ for $i,j\in\hfset$
(\S\ref{sec:weylalgs:Herm}).

\item[$\{a_i^*\}_{i\in\hfset}$]    The dual basis to
$\{a_i\}_{i\in\hfset}$ (\S\ref{sec:weylalgs:Herm}).

\item[$g(\cdot)$]   We write $g\mapsto
g(\cdot)$ for the natural homomorphism $\Mp(\gt{u})\to\Sy(\gt{u})$.
Regarding $g\in\Sy(\gt{u})$ as an element of $\Wy(\gt{u})^{\times}$
we have $g(u)=gug^{-1}$ in $\Wy(\gt{u})$ for $u\in\gt{u}$. More
generally, we write $g(x)$ for $gxg^{-1}$ when $x$ is any element of
$\Wy(\gt{u})$.

\item[$\ogi$]  The map which is $-\Id$ on $\gt{u}$, or the parity
involution on $\Wy(\gt{u})$ (\S\ref{sec:weylalgs:struc}), or the
parity involution on $\WS(\gt{u})_{\Theta}$
(\S\ref{sec:weylalgs:VOAs}).

\item[$\ogih$] The map which is $\ii\Id$ on $\gt{a}$ and $-\ii\Id$
on $\gt{a}^*$, or the lift of this map to $\Wy(\gt{u})$ or
$\WS(\gt{u})_{\Theta}$ (\S\ref{sec:weylalgs:Herm}).

\item[$\hfset$]     A finite ordered set indexing an orthonormal
basis for $\gt{a}$ (\S\ref{sec:weylalgs:Herm}).

\item[$\hfset'$]    A second copy of $\hfset$ with the natural
identification $\hfset\leftrightarrow\hfset'$ denoted
$i\leftrightarrow i'$ for $i\in\hfset$ (\S\ref{sec:weylalgs:Herm}).

\item[$\fset$] A finite ordered set indexing an orthonormal basis
for $\gt{u}$ (\S\ref{sec:weylalgs:struc}). In the case that
$\gt{u}=\gt{a}\oplus\gt{a}^*$, we set $\fset=\hfset\cup\hfset'$ and
insist that $\fset$ be ordered according to the ordering on $\hfset$
and the rule $i<j'$ for $i,j\in\hfset$ (\S\ref{sec:weylalgs:Herm}).

\item[$\mc{E}$]     A label for the basis
$\{e_i\}_{i\in\fset}$ (\S\ref{sec:weylalgs:struc}).

\item[$\WS(\gt{u})$]  The Weyl module VOA associated to the
vector space $\gt{u}$ (\S\ref{sec:weylalgs:VOAs}).

\item[$\WS(\gt{u})_{\ogi}$]     The canonically twisted module
over $A(\gt{u})$ (\S\ref{sec:weylalgs:VOAs}).

\item[$\WS(\gt{u})_{\Theta}$]   The direct sum of $A(\gt{u})$-modules
$A(\gt{u})\oplus A(\gt{u})_{\ogi}$ (\S\ref{sec:weylalgs:VOAs}).

\item[$\Wy(\gt{u})$]     The Weyl algebra associated to the vector space
$\gt{u}$ (\S\ref{sec:weylalgs:struc}).

\item[$\Mp(\gt{u})$]     The metaplectic group associated to the vector
space $\gt{u}$ (\S\ref{sec:weylalgs:meta}).

\item[$\llab\cdot\,,\cdot\rab$]     A non-degenerate antisymmetric
bilinear form on $\gt{u}$ or on $\Wy(\gt{u})$
(\S\ref{sec:weylalgs:struc}).

\item[$(\cdot\,,\cdot)$] A non-degenerate Hermitian form on
$\gt{a}$. The Hermitian forms arising will always be antilinear in
the right hand slot.
\end{list}
\end{small}

\section{Vertex algebras}\label{sec:VOAs}

We briefly review the definition of vertex algebra, and of vertex
operator algebra (VOA) in \S\ref{sec:VOAs:struc}. We recall our
notational conventions regarding vertex Lie algebras in
\S\ref{sec:VOAs:VLAs}. For more background, we refer reader to
\cite[\S1]{DunVARuI} and references therein.

\subsection{Vertex algebra structure}\label{sec:VOAs:struc}

For a {\em vertex algebra} structure on a (super)space
$U=U_{\bar{0}}\oplus U_{\bar{1}}$ over a field $\FF$ we require the
following data.
\begin{itemize}
\item{\em Vertex operators:}   an even morphism $Y:U\otimes U\to
U\ldb z\rdb$ such that when we write
$Y(u,z)v=\sum_{n\in\ZZ}u_{(n)}vz^{-n-1}$, we have $Y(u,z)=0$ only
when $u=0$.
\item{\em Vacuum:}   a distinguished vector $\vac\in U_{\bar{0}}$
such that $Y(\vac,z)u=u$ for $u\in U$, and $Y(u,z)\vac|_{z=0}=u$.
\end{itemize}
This data furnishes a vertex algebra structure on $U$ just when the
following identity is satisfied.
\begin{itemize}
\item{\em Jacobi identity:}   for $\ZZ/2$-homogeneous $u,v\in U$,
and for any $m,n,l\in\ZZ$ we have
\begin{gather}\label{eqn:VOAs:struc:Jac}
     \begin{split}
     &\Res_{z}Y(u,z)Y(v,w)\iota_{z,w}F(z,w)\\
     &-\Res_{z}(-1)^{|u||v|}Y(v,w)Y(u,z)\iota_{w,z}F(z,w)\\
     &=\Res_{z-w}Y(Y(u,z-w)v,w)\iota_{w,z-w}F(z,w)
     \end{split}
\end{gather}
where $F(z,w)=z^mw^n(z-w)^l$.
\end{itemize}
We denote such an object by the triple $(U,Y,\vac)$.

A {\em vertex operator algebra (VOA)} is vertex algebra $(U,Y,\vac)$
equipped with a distinguished element $\cas\in\U_{\bar{0}}$ called
the {\em Virasoro element}, such that $L(-1):=\cas_{0}$ satisfies
\begin{gather}
     [L(-1),Y(u,z)]=D_zY(u,z)
\end{gather}
for all $u\in U$, and such that the operators $L(n):=\cas_{(n+1)}$
furnish a representation of the Virasoro algebra on the vector space
underlying $U$, so that we have
\begin{gather}\label{eqn:VOAs:struc:VirRel}
    \left[L(m),L(n)\right]=(m-n)L(m+n)
    +\frac{m^3-m}{12}\delta_{m+n,0}{ c}{\rm Id}
\end{gather}
for some $c\in\FF$. We also require the following grading condition.
\begin{itemize}
\item{\em $L(0)$-grading:}  the action of $L(0)$ on $U$ is
diagonalizable with rational eigenvalues bounded from below, and is
such that the $L(0)$-homogeneous subspaces $U_n:=\{u\in U\mid
L(0)u=nu\}$ are finite dimensional.
\end{itemize}
When these conditions are satisfied we write $(U,Y,\vac,\cas)$ in
order to indicate the particular data that constitutes the VOA
structure on $U$. The value $c$ in (\ref{eqn:VOAs:struc:VirRel}) is
called the {\em rank} of $U$, and we denote it by $\rank(U)$.

We refer the reader to \cite{FHL} for a discussion of VOA modules,
twisted modules, adjoint operators, and intertwining operators. A
VOA is said to be {\em self-dual} in the case that it has no
non-trivial irreducible modules other than itself.

\subsection{Vertex Lie algebra structure}\label{sec:VOAs:VLAs}

The notion of vertex Lie algebra is an axiomatic formulation of the
kind of object one obtains by replacing the formal series $Y(u,z)v$
with its singular terms $\Sing Y(u,z)v=\sum_{n\geq
0}u_{(n)}vz^{-n-1}$ in the definition of vertex algebra. More
precisely, for a structure of vertex Lie algebra on a superspace
$R=R_{\bar{0}}\oplus R_{\bar{1}}$ we require morphisms $Y_-:R\otimes
R\to z^{-1}R[z^{-1}]$ and $T:R\to R$ such that the following axioms
are satisfied.
\begin{itemize}
\item{\em Translation:} $Y_-(Tu,z)=D_zY(u,z)$.
\item{\em Skew--symmetry:} $Y_-(u,z)v=\Sing e^{zT}Y_-(v,-z)u$.
\item{\em Jacobi identity:} for $\ZZ/2$-homogeneous $u,v\in R$,
and for any $m,n,l\in\ZZ$ we have
\begin{gather}\label{eqn:VOAs:VLAs:Jac}
     \begin{split}
     &\Sing\Res_{z}Y_-(u,z)Y_-(v,w)\iota_{z,w}F(z,w)\\
     &-\Sing\Res_{z}(-1)^{|u||v|}Y_-(v,w)Y_-(u,z)\iota_{w,z}F(z,w)\\
     &=\Sing\Res_{z-w}Y_-(Y_-(u,z-w)v,w)\iota_{w,z-w}F(z,w)
     \end{split}
\end{gather}
where $F(z,w)=z^mw^n(z-w)^l$.
\end{itemize}
We denote such an object by $R=(R,Y_-,T)$. We write the image of
$u\otimes v$ under $Y_-$ as $Y_-(u,z)v=\sum_{n\geq
0}u_{(n)}vz^{-n-1}$. 

For any vertex algebra $U=(U,Y,\vac)$ we obtain a vertex Lie algebra
$(U,Y_-,T)$ by setting $Y_-=\Sing Y$, and by setting $T$ to be the
{\em translation operator:} the morphism on $U$ defined so that
$Tu=u_{-2}\vac$ for $u\in U$. (If $(U,Y,\vac,\cas)$ is a VOA, then
$T$ so defined satisfies $T=L(-1)$.) We abuse notation somewhat to
write $\Sing U$ for this object $(U,\Sing Y,T)$. Conversely, to any
vertex Lie algebra one may canonically associate a vertex algebra
called the {\em enveloping vertex algebra}, and this construction
plays an analogous role for vertex Lie algebras as universal
enveloping algebras do for ordinary Lie algebras.

On the other hand, to each vertex Lie algebra $R$ is canonically
associated a Lie algebra $\Lie(R)$ called the {\em local Lie
algebra} of $R$. As a vector space, we have
$\Lie(R)=R[t,t^{-1}]/\Img \partial$ where $\partial$ is the operator
$T\otimes 1+\Id_R\otimes D_t$ on $R[t,t^{-1}]=
R\otimes\FF[t,t^{-1}]$. Writing $u_{[m]}$ for the image of $u\otimes
t^m$ in $\Lie(R)$ we have
\begin{prop}\label{prop:VOAs:VLAs:loclielag}
$\Lie(R)$ is a Lie algebra under the Lie bracket
\begin{gather}
     \left[u_{[m]},v_{[n]}\right]=\sum_{k\geq 0}\binom{m}{k}
          (u_{(k)}v)_{[m+n-k]}
\end{gather}
and $(Tu)_{[m]}=-mu_{[m-1]}$ for all $u$ in $R$.
\end{prop}
Given a vertex Lie algebra $R$ and a subset $\cgset\subset R$ we may
consider the {\em vertex Lie subalgebra generated by $\cgset$}. This
is by definition just the intersection of all vertex Lie subalgebras
of $R$ that contain $\cgset$. When $\cgset\subset U$ for some vertex
algebra $(U,Y,\vac)$, we will write $[\cgset]$ for the vertex Lie
subalgebra of $\Sing U$ generated by $\cgset$.


\section{Enhanced vertex algebras}\label{sec:eVOAs}

We recall the definitions of enhanced vertex algebra and enhanced
vertex operator algebra in this section. We refer the reader to
\cite[\S2]{DunVARuI} for a fuller discussion.

\begin{defn}
An {\em enhanced vertex algebra} is a quadruple $(U,Y,\vac,R)$ such
that $(U,Y,\vac)$ is a vertex algebra, and $R$ is a vertex Lie
subalgebra of $\Sing(U,Y,\vac)$. We say that $(U,Y,\vac,R)$ is an
{\em enhanced vertex operator algebra (enhanced VOA)} if there is a
unique $\cas\in R$ such that
\begin{enumerate}
\item     $(U,Y,\vac,\cas)$ is a VOA, and
\item     $\cas_{n}u=0$ for all $n\geq 2$ whenever $u\in R$ and
$\cas_{1}u=u$.
\end{enumerate}
When $(U,Y,\vac,R)$ is an enhanced VOA, we call $(U,Y,\vac,\cas)$
the {\em underlying VOA}. The element $\cas$ is called the {\em
Virasoro element} of $(U,Y,\vac,R)$.
\end{defn}
In many instances, the vertex Lie algebra $R$ will be of the form
$R=[\cgset]$ for some subset $\cgset\subset U$, and in such a case
we will write $(U,Y,\vac,\cgset)$ in place of $(U,Y,\vac,[\cgset])$
since there is no loss of information. We will write $(U,\cgset)$ or
even $U$ in place of $(U,Y,\vac,\cgset)$ when no confusion will
arise. When $U=(U,Y,\vac,R)$ is an enhanced vertex algebra we say
that $R$ determines the {\em conformal structure}\footnote{Then
there is a convenient coincidence with the terminology of
\cite{Kac_VAsBegin}, where the objects we refer to as vertex Lie
algebras are called {\em conformal algebras}.} on $U$. The
automorphism group of an enhanced vertex algebra $(U,\cgset)$ is the
subgroup of the group of vertex algebra automorphisms of $U$ that
fixes each element of $\cgset$. In practice, we may write
$\Aut(U,\cgset)$ in order to emphasize this.

For $U=(U,Y,\vac,\cgset)$ an enhanced vertex algebra, the set
$\cgset$ is called a {\em conformal generating set} for $U$, and the
elements of $\cgset$ are called {\em conformal generators}.
\begin{defn} Let $(U,\cgset)=(U,Y,\vac,\cgset)$ be an enhanced
vertex algebra and set $\cfalg(\cgset)=\Lie([\cgset])$, so that
$\cfalg(\cgset)$ is the local Lie algebra of the vertex Lie
subalgebra of $\Sing(U,Y,\vac)$ generated by $\cgset$. We call
$\cfalg(\cgset)$ the {\em local algebra} associated to $(U,\cgset)$.
\end{defn}
The {\em $\U(1)$ Virasoro algebra} is the (purely even) Lie algebra
spanned by symbols $J_m$, $L_m$, and ${\bf{c}}$, for $m\in\ZZ$, and
subject to the following relations.
\begin{gather}\label{eqn:eVOAs:UVirRels}
\begin{split}
    [L_m,L_n]&=(m-n)L_{m+n}+\frac{m^3-m}{12}\delta_{m+n,0}{\bf{c}},
     \quad\left[L_m,{\bf{c}}\right]=0,\\
    [L_m,J_n]&=-nJ_{m+n},\quad
    [J_m,J_n]=-m\delta_{m+n,0}{\bf c},\quad
    \left[J_m,{\bf{c}}\right]=0.
\end{split}
\end{gather}
\begin{defn}
We say that an enhanced VOA $(U,Y,\vac,\cgset)$ is an {\em enhanced
$\U(1)$--VOA} if there is a unique (up to sign) $\cej\in [\cgset]$
such that the Fourier components of the operators $Y(\cas,z)=\sum
L(m)z^{-m-2}$ and $Y(\cej,z)=\sum J(m)z^{-m-1}$ furnish a
representation of the $\U(1)$ Virasoro algebra
(\ref{eqn:eVOAs:UVirRels}) under the assignment $L_m\mapsto L(m)$,
$J_m\mapsto J(m)$. We require also that $J(0)$ acts semi-simply on
$U$.
\end{defn}
Note that for an enhanced $\U(1)$--VOA the element
$\cas_{\alpha}=\cas+\alpha T\cej$ may render
$(U,Y,\vac,\cas_{\alpha})$ a VOA for many choices of $\alpha\in\FF$.
On the other hand, there is only one choice ($\alpha=0$) for which
$(\cas_{\alpha})_{(2)}\cej=0$. Recall from \cite[\S6]{DunVARuI} that
enhanced $\U(1)$--VOAs admit two variable analogues of the usual
McKay--Thompson series that are defined for ordinary VOAs.


\section{Weyl algebras}\label{sec:weylalgs}

The construction of VOAs that we will use arises from certain Weyl
algebra modules. In this section we recall some basic properties of
Weyl algebras. We discuss briefly the metaplectic groups $\Mp_{2N}$
in \S\ref{sec:weylalgs:meta}, and in \S\ref{sec:weylalgs:VOAs} we
recall the construction of VOA module structure on modules over
certain (infinite dimensional) Weyl algebras. In
\S\ref{sec:weylalgs:Herm} we review the Hermitian structure that
arises naturally on these objects given the existence of a suitable
Hermitian form.

\subsection{Weyl algebra structure}\label{sec:weylalgs:struc}

Recall that $\FF$ denotes either $\RR$ or $\CC$. We suppose that
$\gt{u}$ is an $\FF$-vector space of even dimension with
non-degenerate antisymmetric bilinear form $\llab
\cdot\,,\cdot\rab$.

We write $\Wy(\gt{u})$ for the Weyl algebra over $\FF$ generated by
$\gt{u}$. More precisely, we set $\Wy(\gt{u})=T(\gt{u})/I(\gt{u})$
where $T(\gt{u})$ is the tensor algebra of $\gt{u}$ over $\FF$ with
unit denoted ${\bf 1}$, and $I(\gt{u})$ is the ideal of $T(\gt{u})$
generated by all expressions of the form $u\otimes v-v\otimes
u+2\llab u,v\rab{\bf 1}$ for $u,v\in \gt{u}$. The natural algebra
structure on $T(\gt{u})$ induces an associative algebra structure on
$\Wy(\gt{u})$. The vector space $\gt{u}$ embeds naturally in
$\Wy(\gt{u})$, and when it is convenient we identify $\gt{u}$ with
its image in $\Wy(\gt{u})$. We also write $\alpha$ in place of
$\alpha{\bf 1}+I(\gt{u})\in\Wy(\gt{u})$ for $\alpha\in\FF$ when no
confusion will arise. For $u,v\in \gt{u}$ we have the relation
$uv-vu=-2\llab u,v\rab$ in $\Wy(\gt{u})$, and more generally
\begin{gather}
     u^mv^n=\sum_{k\geq 0}\binom{m}{k}\binom{n}{k}
          2^k\llab v,u\rab^k v^{n-k}u^{m-k}
\end{gather}
for nonnegative integers $m$ and $n$.

The linear transformation on $\gt{u}$ which is $-1$ times the
identity map lifts naturally to $T(\gt{u})$ and preserves
$I(\gt{u})$, and hence induces an involution on $\Wy(\gt{u})$ which
we denote by $\ogi$. The map $\ogi$ is known as the {\em parity
involution}. We have $\ogi(u_1\cdots u_k)=(-1)^ku_1\cdots u_k$ for
$u_1\cdots u_k\in\Wy(\gt{u})$ with $u_i\in \gt{u}$, and we write
$\Wy(\gt{u})=\Wy(\gt{u})^0\oplus \Wy(\gt{u})^1$ for the
decomposition into eigenspaces for $\ogi$.

A linear isomorphism $\lambda:\bigvee(\gt{u})\to\Wy(\gt{u})$ is
defined by setting
\begin{gather}
     \lambda(u_1\vee\cdots\vee u_k)
     =\frac{1}{k!}\sum_{\sigma\in S_k}
          u_{\sigma 1}\cdots u_{\sigma k}\vac
\end{gather}
for $u=u_1\vee\cdots\vee u_k$ with $u_i\in \gt{u}$. The
antisymmetric bilinear form on $\gt{u}$ extends naturally to the
symmetric algebra $\bigvee(\gt{u})$ by decreeing that
$\bigvee^m(\gt{u})$ be orthogonal to $\bigvee^n(\gt{u})$ for $m\neq
n$, and by setting
\begin{gather}
     \llab u_1\vee\cdots\vee u_k,v_1\vee\cdots\vee v_k\rab
     =\sum_{\sigma\in S_k}
          \llab u_{\sigma 1},v_1\rab\cdots
          \llab u_{\sigma k},v_k\rab
\end{gather}
and $\llab \vac,\vac\rab=1$, for $\vac$ the unit in
$\bigvee(\gt{u})$. Evidently, this form is antisymmetric on odd
symmetric powers of $\gt{u}$, and symmetric on even powers. Now we
may use the map $\lambda$ to imbue $\Wy(\gt{u})$ with a bilinear
form by setting
\begin{gather}
     \llab a,b\rab=\llab \lambda^{-1}(a),\lambda^{-1}(b)\rab
\end{gather}
for $a,b\in\Wy(\gt{u})$. With this definition we have
$\llab\vac,\vac\rab=1$ and $\llab ua,b\rab=-\llab a,ub\rab$ for
$a,b\in \Wy(\gt{u})$ and $u\in\gt{u}$, and the restriction to
$\gt{u}\hookrightarrow\Wy(\gt{u})$ agrees with the original form on
that space.

\begin{eg}
We can compute $\llab\vac,uv\rab$ for $u,v\in\gt{u}$, by observing
that
\begin{gather}
     \lambda(u\vee v)
          =\frac{1}{2}(uv+vu)
          =uv+\llab u,v\rab
\end{gather}
in $\Wy(\gt{u})$, so that $\llab \vac,u\vee v\rab=0$ in
$\bigvee(\gt{u})$ implies $\llab\vac,uv\rab=-\llab u,v\rab$.
\end{eg}

\subsection{Metaplectic groups}\label{sec:weylalgs:meta}

We would like to define an enlargement of $\Wy(\gt{u})$ that will
contain expressions like $\exp(tu^2)$ for $t\in\FF$ and
$u\in\gt{u}$. Our construction follows \cite{Cru_OrthSympCliffAlgs}.

Let $\{e_i,f_{i}\}_{i=1}^N$ be a basis for $\gt{u}$ satisfying
$\llab e_i,f_j\rab=\delta_{ij}$. For $H=(h_1,\ldots,h_N)$ in
$\ZZ_{\geq 0}^N$ we write $e^H$ for $e_1^{h_1}\cdots e_N^{h_N}$, and
similarly for $f^H$. Then $\Wy(\gt{u})$ is spanned by the
expressions $e^Hf^K$ for $H,K\in\ZZ_{\geq 0}^N$. We set $H!=\prod
h_i!$ and $|H|=\sum h_i$. We define $\Wyx(\gt{u})$ to be the set of
expressions of the form
\begin{gather}
     \widetilde{u}=\sum_{H,K}\frac{M_{H,K}e^Hf^K}
          {(H!K!)^{1/2}}
\end{gather}
where the $M_{H,K}\in\FF$ should satisfy $|M_{H,K}|\leq C/
2^{|H|+|K|}$ for some constant $C$, for all $H$ and $K$ such that
$|H|$ and $|K|$ are sufficiently large. Then $\Wy(\gt{u})$ embeds
naturally in $\Wyx(\gt{u})$, and the latter space is closed under
the algebra structure it naturally inherits from the former.

Let $\Wyx(\gt{u})^{\times}$ denote the set of invertible elements in
$\Wyx(\gt{u})$. For $x\in\Wyx(\gt{u})^{\times}$ and
$a\in\Wy(\gt{u})$ we set $x(a)=xax^{-1}$. We define the {\em
metaplectic group associated to $\gt{u}$} to be the subgroup
$\Mp(\gt{u})$ of $\Wyx(\gt{u})$ generated by expressions of the form
$\alpha\exp(tu^2)$ for $\alpha,t\in\FF$ and $u\in\gt{u}$, with
$|\alpha|=1$. The identity
\begin{gather}
     \exp(tu^2)v\exp(-tu^2)=v-4\llab u,v\rab tu
\end{gather}
for $v\in \gt{u}$ shows that $x(\cdot)$ preserves $\gt{u}$ for
$x\in\Mp(\gt{u})$, and further, that the map $x\mapsto x(\cdot)$
defines a surjection $\Mp(\gt{u})\to\Sy(\gt{u})$. We have a short
exact sequence
\begin{gather}
     1\to T\to\Mp(\gt{u})\to\Sy(\gt{u})\to 1
\end{gather}
where $T$ denotes the elements of unit norm in $\FF$.

\subsection{Weyl module VOAs}\label{sec:weylalgs:VOAs}

In this section we review the construction of VOA structure on
modules over a certain infinite dimensional Weyl algebra associated
to $\gt{u}$. The construction appears also in \cite{WeiSympPVOA},
and we refer the reader there for further details. Another important
precursor is the article \cite{FeiFreClsAffAlgs} in which the
construction below, and generalizations thereof, are applied to the
theory of affine Lie (super)algebras.

\medskip

Suppose that $\gt{u}$ has dimension $2N$. Let $\hfset$ be a set of
cardinality $N$, and let $\hfset'$ be a second copy of this set,
with the natural correspondence denoted $i\leftrightarrow i'$. (The
inverse of the map $i\mapsto i'$ will also be denoted by $'$, so
that $i''=i$.) We assume chosen a basis
$\mc{E}=\{e_i,e_{i'}\}_{i\in\hfset}$ for $\gt{u}$ such that
\begin{gather}
     \llab e_i,e_{j'}\rab=-\llab e_{i'}e_j\rab=\delta_{ij}
\end{gather}
Let $\hat{\gt{u}}$ and $\hat{\gt{u}}_{\ogi}$ denote the infinite
dimensional vector spaces with antisymmetric bilinear form described
as follows.
\begin{gather}
    \hat{\gt{u}}=\coprod_{m\in\ZZ}\gt{u}\otimes t^{m+1/2},\quad
    \hat{\gt{u}}_{\ogi}=\coprod_{m\in\ZZ}\gt{u}\otimes t^m,\\
    \llab u\otimes t^r,v\otimes t^s\rab
        =\llab u,v\rab \delta_{r+s,0},
    \; \text{ for $u,v\in\gt{u}$ and $r,s\in\hZZ$.}
\end{gather}
We write $u(r)$ for $u\otimes t^r$ when $u\in\gt{u}$ and $r\in\hZZ$.
We consider the Weyl algebras $\Wy(\hat{\gt{u}})$ and
$\Wy(\hat{\gt{u}}_{\ogi})$. By the conventions of
\S\ref{sec:weylalgs:struc} we have the following {\em fundamental
relation}
\begin{gather}\label{eqn:weylalgs:VOAs:fundid}
     [u(r),u'(s)]=-2\llab u,u'\rab\delta_{r+s,0}
\end{gather}
for any $u,u'\in\gt{u}$, holding in $\Wy(\hat{\gt{u}})$ for
$r,s\in\ZZh$, and holding in $\Wy(\hat{\gt{u}}_{\ogi})$ for
$r,s\in\ZZ$.
\begin{rmk}
Our generators $e_i(r)$ and $e_{i'}(r)$ correspond to
$\sqrt{2}a_i^-(r)$ and $\sqrt{2}a_i^+(r)$, respectively, in the
notation of \cite{WeiSympPVOA}.
\end{rmk}
The inclusion of $\gt{u}$ in $\hat{\gt{u}}_{\ogi}$ given by
$u\mapsto u(0)$ induces an embedding of algebras
$\Wy(\gt{u})\hookrightarrow \Wy(\hat{\gt{u}}_{\ogi})$.

We write $\mc{B}(\hat{\gt{u}})$ for the subalgebra of
$\Wy(\hat{\gt{u}})$ generated by the $u(m+\tfrac{1}{2})$ for
$u\in\gt{u}$ and $m\in\ZZ_{\geq 0}$. We write
$\mc{B}(\hat{\gt{u}}_{\ogi})$ for the subalgebra of
$\Wy(\hat{\gt{u}}_{\ogi})$ generated by the $u(m)$ for $u\in\gt{u}$
and $m\in\ZZ_{> 0}$, and by the $u(0)$ for $u\in\gt{a}^*$. Let
$\FF_0$ denote a one-dimensional module for either
$\mc{B}(\hat{\gt{u}})$ or $\mc{B}(\hat{\gt{u}}_{\ogi})$, spanned by
a vector $1_0$, such that $u(r)1_0=0$ whenever $r\in\hZZ_{>0}$, and
such that $u(0)1_{0}={0}$ for $u\in \gt{a}^*$. We write
$\WS(\gt{u})$ (respectively $\WS(\gt{u})_{\ogi}$) for the
$\Wy(\hat{\gt{u}})$-module (respectively
$\Wy(\hat{\gt{u}}_{\ogi})$-module) induced from the
$\mc{B}(\hat{\gt{u}})$-module structure (respectively
$\mc{B}(\hat{\gt{u}}_{\ogi})$-module structure) on $\FF_{0}$.
\begin{gather}
    \WS(\gt{u})=\Wy(\hat{\gt{u}})
        \otimes_{\mc{B}(\hat{\gt{u}})}\FF_{0},\qquad
    \WS(\gt{u})_{\ogi}=\Wy(\hat{\gt{u}}_{\ogi})
        \otimes_{\mc{B}(\hat{\gt{u}}_{\ogi})}\FF_{0}.
\end{gather}
We write $\vac$ for the vector $1\otimes 1_{0}$ in $\WS(\gt{u})$,
and we write $\vac_{\ogi}$ for the vector $1\otimes 1_{0}$ in
$\WS(\gt{u})_{\ogi}$.

\medskip

The space $\WS(\gt{u})$ supports a structure of VOA. In order to
define the vertex operators we review the notion of bosonic normal
ordering for elements in $\Wy(\hat{\gt{u}})$ and
$\Wy(\hat{\gt{u}}_{\ogi})$. The bosonic normal ordering on
$\Wy(\hat{\gt{u}})$ is the multi-linear operator defined so that for
$u_i\in\gt{u}$ and $r_i\in\ZZh$ we have
\begin{equation}
    :\!u_1(r_1)\cdots u_k(r_k)\!:\;=\,
        u_{\sigma 1}(r_{\sigma 1})
        \cdots u_{\sigma k}(r_{\sigma k})
\end{equation}
where $\sigma$ is any permutation of the index set $\{1,\ldots,k\}$
such that $r_{\sigma 1}\leq\cdots\leq r_{\sigma k}$. For elements in
$\Wy(\hat{\gt{u}}_{\ogi})$ the bosonic normal ordering is defined in
steps by first setting
\begin{equation}
    :\!u_1(0)\cdots u_k(0)\!:\;=\frac{1}{k!}\sum_{\sigma\in S_k}
        u_{\sigma 1}(0)
        \cdots u_{\sigma k}(0)
\end{equation}
for $u_i\in \gt{u}$. Then in the situation that $n_i\in\ZZ$ are such
that $n_{i}\leq n_{i+1}$ for all $i$, and there are some $s$ and $t$
(with $1\leq s\leq t\leq k$) such that $n_j=0$ for $s\leq j\leq t$,
we set
\begin{equation}
     \begin{split}
    &:\!u_1(n_1)\cdots u_k(n_k)\!:\\
    =u_{1}(n_{1})\cdots u_{ s-1}(n_{s-1})&
    \!:\!u_s(0)\cdots u_t(0)\!:\!
    u_{t+1}(n_{t+1})\cdots u_{k}(n_k)
    \end{split}
\end{equation}
Finally, for arbitrary $n_i\in\ZZ$ we set
\begin{equation}
    :\!u_1(n_1)\cdots u_k(n_k)\!:\,=\,
        :\!u_{\sigma 1}(n_{\sigma 1})
        \cdots u_{\sigma k}(n_{\sigma k})\!:
\end{equation}
where $\sigma$ is again any permutation of the index set
$\{1,\ldots,k\}$ such that $n_{\sigma 1}\leq\cdots\leq n_{\sigma
k}$, and we extend the definition multi-linearly to
$\Wy(\hat{\gt{u}}_{\ogi})$.

\medskip

For $u\in\gt{u}$ we now define the generating function, denoted
$u(z)$, of operators on $\WS(\gt{u})_{\Theta}=\WS(\gt{u})\oplus
\WS(\gt{u})_{\ogi}$ by setting
\begin{gather}
     u(z)=\sum_{r\in\hZZ}u(r)z^{-r-1/2}
\end{gather}
Note that $u(r)$ acts as $0$ on $\WS(\gt{u})$ if $r\in \ZZ$, and
acts as $0$ on $\WS(\gt{u})_{\ogi}$ if $r\in \ZZh$. To an element
$a\in \WS(\gt{u})$ of the form $a=u_{1}(-m_1-\tfrac{1}{2})\cdots
u_{k}(-m_k-\tfrac{1}{2}){\bf 1}$ for $u_i\in \gt{u}$ and $m_i\in
\ZZ_{\geq 0}$, we associate the operator valued power series
$\overline{Y}(a,z)$, given by
\begin{gather}
    \overline{Y}(a,z)=\,:\!D_z^{(m_1)}u_{i_1}(z)\cdots
        D_z^{(m_k)}u_{i_k}(z)\!:
\end{gather}
We define the vertex operator correspondence
\begin{gather}
    Y(\cdot\,,z):\WS(\gt{u})\otimes \WS(\gt{u})_{\Theta}
        \to \WS(\gt{u})_{\Theta}((z^{1/2}))
\end{gather}
by setting $Y(a,z)b=\overline{Y}(a,z)b$ when $b\in \WS(\gt{u})$, and
by setting $Y(a,z)b=\overline{Y}(e^{\Delta_z}a,z)b$ when $b\in
\WS(\gt{u})_{\ogi}$, where $\Delta_z$ is the expression defined by
\begin{gather}
    \Delta_z=-\frac{1}{2}\sum_i\sum_{m,n\in\ZZ_{\geq 0}}C_{mn}
        e_{i'}(m+\tfrac{1}{2})e_i(n+\tfrac{1}{2})z^{-m-n-1}\\
        C_{mn}=\frac{1}{2}\frac{(m-n)}{m+n+1}
            \binom{-\tfrac{1}{2}}{m}\binom{-\tfrac{1}{2}}{n}
\end{gather}
Define $\cas$ by setting
\begin{gather}
     \cas=
     \frac{1}{4}\sum_i\left(
     e_{i'}(-\tfrac{1}{2})e_i(-\tfrac{3}{2})\vac
     -e_i(-\tfrac{1}{2})e_{i'}(-\tfrac{3}{2})\vac
     \right)
     \in \WS(\gt{u})_2
\end{gather}
Then we have the following
\begin{thm}[\cite{WeiSympPVOA}]\label{thm:weylalgs:VOAs:VOAstruc}
The map $Y$ defines a structure of VOA of rank $-N$ on $\WS(\gt{u})$
when restricted to $\WS(\gt{u})\otimes \WS(\gt{u})$, and the
Virasoro element is given by $\cas$. The map $Y$ defines a structure
of $\ogi$-twisted $\WS(\gt{u})$-module on $\WS(\gt{u})_{\ogi}$ when
restricted to $\WS(\gt{u})\otimes \WS(\gt{u})_{\ogi}$.
\end{thm}
Given $I=(i_1,\cdots,i_N)\in\fset^N$ with $\fset=\hfset\cup\hfset'$,
let us write $e_I(-\tfrac{1}{2})$ for $e_{i_1}(-\tfrac{1}{2})\cdots
e_{i_N}(-\tfrac{1}{2})$. Observe then that $\WS(\gt{u})_2$ is
spanned by vectors of the form $e_I(-\tfrac{1}{2})\vac$ for
$I\in\fset^4$, and by the $e_i(-\tfrac{3}{2})e_j(-\tfrac{1}{2})\vac$
with $i,j\in\fset$.

All we need to know about the expressions $Y(a,z)b$ for $b\in
\WS(\gt{u})_{\ogi}$ is contained in the following
\begin{prop}\label{prop:weylalgs:VOAs:twops}
Let $b\in \WS(\gt{u})_{\ogi}$.
\begin{enumerate}
\item   If $a=\vac\in \WS(\gt{u})_0$ then $Y(a,z)b=b$.
\item   If $a\in \WS(\gt{u})_1$ then $\Delta_za=0$ so that
$Y(a,z)b=\overline{Y}(a,z)b$.
\item   If $a\in \WS(\gt{u})_2$ and
     $a\in\Span\left\{e_I(-\tfrac{1}{2})\vac,
          e_i(-\tfrac{3}{2})e_j(-\tfrac{1}{2})\vac\mid i'\neq j\right\}$
then $\Delta_za=0$ and $Y(a,z)b= \overline{Y}(a,z)b$.
\item   For $a=e_{i'}(-\tfrac{1}{2})e_i(-\tfrac{3}{2})\vac$ and $i\in\hfset$ we have
$\Delta_za=-\tfrac{1}{4}z^{-2}$ and $\Delta_z^2a=0$ so that
$Y(a,z)b=\overline{Y}(a,z)b-\tfrac{1}{4}bz^{-2}$ in this case.
\item   For $a=e_i(-\tfrac{1}{2})e_{i'}(-\tfrac{3}{2})\vac$ and $i\in\hfset$ we have
$\Delta_za=\tfrac{1}{4}z^{-2}$ and $\Delta_z^2a=0$ so that
$Y(a,z)b=\overline{Y}(a,z)b+\tfrac{1}{4}bz^{-2}$ in this case.
\end{enumerate}
\end{prop}
As a corollary of Proposition~\ref{prop:weylalgs:VOAs:twops} we have
that
\begin{gather}
     Y(\omega,z){\bf 1}_{\ogi}=-\frac{N}{8}
          {\bf 1}_{\ogi} z^{-2}
\end{gather}
and consequently the $L(0)$-grading on $\WS(\gt{u})_{\Theta}$ is
given as follows.
\begin{gather}
    \WS(\gt{u})=\coprod_{n\in\hZZ_{\geq 0}}\WS(\gt{u})_n,\quad
    \WS(\gt{u})_{\ogi}=\coprod_{n\in\hZZ_{\geq 0}}
        (\WS(\gt{u})_{\ogi})_{n-N/8}.
\end{gather}
Note that $\WS(\gt{u})$ is {\em strongly generated} (c.f.
\cite{Kac_VAsBegin}) by it's degree $1/2$ subspace in the sense that
we have
\begin{gather}\label{eqn:weylalgs:VOAs:strnggens}
     \WS(\gt{u})=\Span\left\{u^1(-m_1-\tfrac{1}{2})\cdots
     u^k(-m_k-\tfrac{1}{2})\vac\mid
          u^i\in \gt{u},\;m_i\in\ZZ_{\geq 0}\right\}
\end{gather}
Armed with this observation it is not difficult to see that
$\WS(\gt{u})$ is a self-dual VOA.
\begin{prop}\label{prop:weylalgs:VOAs:selfdual}
Let $M$ be an irreducible $\WS(\gt{u})$-module with grading bounded
from below. Then $M$ is isomorphic to $\WS(\gt{u})$ as a
$\WS(\gt{u})$-module. In particular, $\WS(\gt{u})$ is a self-dual
VOA.
\end{prop}
\begin{proof}
Let $v_0\in M$ be a homogeneous non-zero element of minimal degree.
Then we have $u(m+\tfrac{1}{2})v_0=0$ for any $u\in\gt{u}$ whenever
$m\geq 0$. We claim that $u(-m-\tfrac{1}{2})v_0\neq 0$ for any
non-zero $u\in \gt{u}$ and for all $m\geq 0$ since otherwise we have
\begin{gather}
     0=u'(m+\tfrac{1}{2})u(-m-\tfrac{1}{2})v_0=
          [u'(m+\tfrac{1}{2}),u(-m-\tfrac{1}{2})]v_0=\llab u,u'\rab v_0
\end{gather}
with the second equality holding because $u'(m+\tfrac{1}{2})v_0=0$.
Evidently then, the map $\WS(\gt{u})\to M$ given by
\begin{gather}
     u^1(-m_1-\tfrac{1}{2})\cdots u^k(-m_k-\tfrac{1}{2})\vac
     \mapsto
     u^1(-m_1-\tfrac{1}{2})\cdots u^k(-m_k-\tfrac{1}{2})v_0
\end{gather}
for $u^i\in\gt{u}$ and $m_i\in\ZZ_{\geq 0}$ is an embedding of
$\WS(\gt{u})$-modules. (The map is defined on all of $\WS(\gt{u})$
by (\ref{eqn:weylalgs:VOAs:strnggens}).) Since $M$ is irreducible by
hypothesis, the image of this map is all of $M$, and the proposition
follows.
\end{proof}

The metaplectic group $\Mp(\gt{u})$ acts naturally on
$\WS(\gt{u})_{\Theta}$. This action is generated by exponentials of
the operators $x_{0}$ for $x\in \WS(\gt{u})_1$. In particular, any
$a\in\Mp(\gt{u})\subset\Wyx(\gt{u})$ may be regarded as a VOA
automorphism of $\WS(\gt{u})$, and as an equivariant linear
isomorphism of the $\WS(\gt{u})$-module $\WS(\gt{u})_{\ogi}$.

\subsection{Hermitian structure}\label{sec:weylalgs:Herm}

In this section we take special interest in the case that $\FF=\CC$
and $\gt{u}$ is of the form $\gt{u}=\gt{a}\oplus\gt{a}^*$ for
$\gt{a}$ a complex vector space with non-degenerate Hermitian form,
and $\gt{a}^*$ the dual space to $\gt{a}$. The Hermitian form will
be denoted $(\cdot\,,\cdot)$, and we assume it to be antilinear in
the second slot. We extend $(\cdot\,,\cdot)$ in the natural way to a
Hermitian form on $\gt{u}$. Suppose that $\gt{a}$ has dimension $N$.
We take $\hfset$ to be some set of cardinality $N$, and we use it to
index some orthonormal basis $\{a_i\}_{i\in\hfset}$ for $\gt{a}$. We
then let $\{a_i^*\}_{i\in\hfset}$ denote the dual basis, so that the
expressions
\begin{gather}
     (\alpha a_i,\beta a_j) =(\alpha a_i^*,\beta a_j^*)
     =\alpha\bar{\beta}\delta_{ij},\quad
     (a_i,a_j^*)=(a_i^*,a_j)=0,
\end{gather}
completely describe the Hermitian form on $\gt{u}$.

There is a natural symmetric bilinear form on $\gt{u}$ induced from
the pairing between $\gt{a}$ and $\gt{a}^*$. We take
$\lab\cdot\,,\cdot\rab$ to be $1/2$ times this form. More precisely,
we set $\lab a,f\rab =\lab f,a\rab =\tfrac{1}{2}f(a)$ for
$a\in\gt{a}$ and $f\in\gt{a}^*$. We make the convention that
$\mc{E}$ shall denote the basis for $\gt{u}$ obtained in the
following way. We let $\hfset'$ be another copy of the set $\hfset$
with the natural identification $\hfset\leftrightarrow\hfset'$
denoted $i\leftrightarrow i'$, and we set
$\mc{E}=\{e_i,e_{i'}\}_{i\in\hfset}$ where
\begin{gather}
     e_i=a_i+a_i^*,\quad e_{i'}=\ii(a_i-a_i^*),
\end{gather}
for $i\in\hfset$. One can check that the basis $\mc{E}$ is
orthonormal with respect to the symmetric form
$\lab\cdot\,,\cdot\rab$.

Let $\ogih$ be the map which is multiplication by $\ii$ on $\gt{a}$,
and multiplication by $-\ii$ on $\gt{a}^*$. Then $\ogih$ preserves
both the Hermitian form $(\cdot\,,\cdot)$ and the symmetric form
$\lab\cdot\,,\cdot\rab$, and we may define an antisymmetric form
$\llab\cdot\,,\cdot\rab$ on $\gt{u}$ by setting
\begin{gather}
     \llab u,v\rab=\lab\ogih u,v\rab
\end{gather}
for $u,v\in\gt{u}$. (Mnemonic: ``$\llab$'' is a shorthand for
``$\lab\ogih$''.) We then have $\llab e_i,e_{j'}\rab=-\llab
e_{i'},e_j\rab=\delta_{ij}$ for $i,j\in\hfset$. (Compare with the
convention of \S\ref{sec:weylalgs:VOAs}.) We also have $\llab
a_i,a_{j}^*\rab=-\llab a_{i}^*,a_j\rab=\tfrac{1}{2}\ii\delta_{ij}$
so that
\begin{gather}
     a_i^*a_j=\ii\delta_{ij}+a_ja_i^*
\end{gather}
in $\Wy(\gt{u})$, by the relation (\ref{eqn:weylalgs:VOAs:fundid}).


\section{Linear groups}\label{sec:GLn}

In this section we define a family of enhanced VOAs with symmetry
groups of the form $\GL_N(\CC)$ for $N$ a positive integer. In
defining these objects we are laying some of the ground work for the
construction of an enhanced VOA for the Rudvalis group in
\S\ref{sec:Ru}.

\medskip

Let $N$ be a positive integer and let $\gt{a}$ be a complex vector
space of dimension $N$ equipped with a positive definite Hermitian
form denoted $(\cdot\,,\cdot)$. As in \S\ref{sec:weylalgs:Herm} we
set ${\gt{u}}=\gt{a}\oplus\gt{a}^*$, we extend the Hermitian form to
$\gt{u}$, and equip this space also with the antisymmetric bilinear
form $\llab\cdot\,,\cdot\rab$, naturally determined as in
\S\ref{sec:weylalgs:Herm}.

Let us denote $\wcgU=\{\cas,\cej\}$ where $\cas$ and $\cej$ are
given by
\begin{gather}
     \begin{split}
     \cas&
          =\frac{1}{4}\sum_{i\in\hfset}
     \left(e_{i'}(-\tfrac{1}{2})e_i(-\tfrac{3}{2})\vac
     -e_{i}(-\tfrac{1}{2})e_{i'}(-\tfrac{3}{2})
          \vac\right)\\
          &
          =\frac{\ii}{2}\sum_{i\in\hfset}
     \left(a_i(-\tfrac{1}{2})a^*_i(-\tfrac{3}{2})\vac
     -a_i^*(-\tfrac{1}{2})a_i(-\tfrac{3}{2})\vac\right)
          \in \WS(\gt{u})_2
     \end{split}\\
     \cej
          =
          \sum_{i\in\hfset}
     a_i(-\tfrac{1}{2})a^*_{i}(-\tfrac{1}{2})\vac\in \WS(\gt{u})_1,
\end{gather}
so that $\cas$ is just the usual Virasoro element for $\WS(\gt{u})$,
and $\cej$ is an element dependent upon the Hermitian structure on
$\gt{a}$. We define operators $L(n)$ and $\ceJ(n)$ for $n\in\ZZ$ by
setting
\begin{align}
     Y(\cas,z)&=\;L(z)=\sum L(n)z^{-n-2},\\
     Y(\cej,z)&=\;\ceJ(z)=\sum \ceJ(n)z^{-n-1}.
\end{align}
We then have
\begin{prop}\label{prop:GLn:GLnOPEs}
The operators $L(z)$ and $\ceJ(z)$ satisfy the following OPEs.
\begin{gather}
     L(z)L(w)=
          \frac{-N/2}{(z-w)^4}+\frac{2L(w)}{(z-w)^2}
          +\frac{D_wL(w)}{(z-w)}
          +{\rm reg.}\\
     L(z)J(w)=\frac{J(w)}{(z-w)^2}+\frac{D_wJ(w)}{(z-w)}
          +{\rm reg.}\\
     J(z)J(w)=\frac{N}{(z-w)^2}
          +{\rm reg.}
\end{gather}
\end{prop}
Proposition~\ref{prop:GLn:GLnOPEs} encodes the commutation relations
amongst the operators $L(m)$ and $J(m)$. We record these in the next
\begin{cor}\label{cor:GLn:LJcomrels}
The operators $L(m)$ and $J(m)$ satisfy the following commutation
relations for $m,n\in\ZZ$.
\begin{gather}
     \left[L(m),L(n)\right]
          =(m-n)L(m+n)
          -\frac{m^3-m}{12}{N}\delta_{m+n,0}\Id\\
     [L(m),J(n)]
          =-nJ(m+n)\\
     [J(m),J(n)]
          =
               mN\delta_{m+n,0}\Id
\end{gather}
\end{cor}
Proposition~\ref{prop:GLn:GLnOPEs} shows that the set
$\wcgU=\{\cas,\cej\}$ furnishes an enhanced conformal structure on
$\WS(\gt{u})$, with respect to which $\wcgU$ is a generating set of
$0$ defect (c.f. \cite[\S2]{DunVARuI}). Note also that the VOA
$\WS(\gt{u})$ has negative rank $-N$, where $N=\dim(\gt{a})$.

\begin{prop}\label{prop:GLn:u(1)struc}
For $N$ a positive integer, the quadruple $(\wu,Y,\vac,\wcgU)$ is a
self-dual $\U(1)$--VOA of rank $-N$.
\end{prop}
\begin{proof}
We will show that $(\wu,\wcgU)$ admits a $\U(1)$--VOA structure
(c.f.~\S\ref{sec:eVOAs}). The fact that the $4$-tuple
$(\wu,Y,\vac,\cas)$ is a VOA of rank $-N$ is the content of
Theorem~\ref{thm:weylalgs:VOAs:VOAstruc}, and self-duality was shown
in Proposition \ref{prop:weylalgs:VOAs:selfdual}.

From Proposition~\ref{prop:GLn:GLnOPEs} we see that $\wcgU$ has
defect $0$, so that the vertex Lie algebra generated by $\wcgU$ is
just $[\wcgU]=\Span\{\vac,T^k\cej,T^k\cas\mid k\in\ZZ\}$ (c.f.
\cite[\S2]{DunVARuI}). From this explicit spanning set we see that
the only Virasoro elements for $\WS(\gt{u})$ in $[\wcgU]$ are of the
form $\cas_{\alpha}=\cas+\alpha T\cej$ for some $\alpha\in \CC$.
(The operator $T$ maps $u=u_{(-1)}\vac$ to $u_{(-2)}\vac$.) Using
Corollary~\ref{cor:GLn:LJcomrels} we check that the Fourier
coefficients of $Y(\cas_{\alpha},z)=L_{\alpha}(m)z^{-m-2}$ furnish a
representation of the Virasoro algebra of rank $-N(1+12\alpha^2)$,
and we have $L_{\alpha}(n)=L(n)+\alpha(-n-1)J(n)$ where
$Y(\cas,z)=\sum L(m)z^{-m-2}$ and $Y(\cej,z)=\sum J(m)z^{-m-1}$. In
particular, $L_{\alpha}(-1)=L(-1)=T$ for all $\alpha$, and
$L_{\alpha}(0)=L(0)-\alpha J(0)$. Thus if $u\in [\wcgU]$ and
$L_{\alpha}(0)u=u$ then $u\in\Span\{\cej\}$. Now we compute
$L_{\alpha}(1)\cej= \alpha(-2)J(1)\cej= 2\alpha N$, and it follows
that $(\wu,Y,\vac,\wcgU)$ is an enhanced VOA with the Virasoro
element given by $\cas=\cas_{0}$. From Corollary
~\ref{cor:GLn:LJcomrels} we see that $\cej$ is the unique (up to
sign) choice of element in $[\cgU]$ for which the component
operators of $Y(\cej,z)$ satisfy the required relations
(\ref{eqn:eVOAs:UVirRels}), and thus $(\WS(\gt{u}),\wcgU)$ is a
$\U(1)$--VOA.
\end{proof}
We would like to compute the automorphism group of
$(\WS(\gt{u}),\wcgU)$. First we will compute the automorphism group
of the underlying VOA $(\WS(\gt{u}),\cas)$.
\begin{prop}\label{prop:GLn:VOAauts}
The automorphism group of $(\WS(\gt{u}),\cas)$ is $\Sy(\gt{u})$.
\end{prop}
\begin{proof}
Observe that there is a natural correspondence
$\gt{u}\leftrightarrow \WS(\gt{u})_{1/2}$ given by $u\leftrightarrow
u(-\tfrac{1}{2})\vac$, and that we recover the corresponding
antisymmetric form on $\WS(\gt{u})_{1/2}$ by defining $\llab
a|b\rab$ for $a,b\in\WS(\gt{u})_{1/2}$ in such a way that
$a_{0}b=\llab a|b\rab\vac$. (The skew-symmetry of the form follows
from the skew-symmetry property of a VOA:
$Y(u,z)v=e^{zL(-1)}Y(v,-z)u$.) Any automorphism of the VOA
$\WS(\gt{u})$ must preserve the homogeneous subspace
$\WS(\gt{u})_{1/2}$, and must preserve the bracket $[u,v]=u_{0}v$,
and thus we obtain a map
$\phi:\Aut(\WS(\gt{u}),\cas)\to\Sy(\gt{u})$. On the other hand,
$\WS(\gt{u})$ is strongly generated by its subspace of degree $1/2$
(c.f. (\ref{eqn:weylalgs:VOAs:strnggens})), in the sense that we
have
\begin{gather}
     \WS(\gt{u})=\Span\left\{u^{1}_{-n_1}\cdots u^k_{-n_k}\vac\mid
          u^i\in \WS(\gt{u})_{1/2},\,n_i> 0,\, k\geq 0\right\}
\end{gather}
and we obtain a map $\psi:\Sy(\gt{u})\to\Aut(\WS(\gt{u}),\cas)$, by
letting $g\in\Sy(\gt{u})$ act on $\WS(\gt{u})$ in the following way.
\begin{gather}
     g:u^{1}_{-n_1}\cdots u^k_{-n_k}\vac\mapsto
          (gu^{1})_{-n_1}\cdots (gu^k)_{-n_k}\vac
\end{gather}
Evidently, $\psi$ is an inverse to $\phi$.
\end{proof}

\begin{prop}\label{prop:GLn:eVOAauts}
The automorphism group of $(\WS(\gt{u}),\wcgU)$ is $\GL(\gt{a})$.
\end{prop}
\begin{proof}
Recall that $\cej= \sum a_i^*(-\tfrac{1}{2})a_i(-\tfrac{1}{2})\vac$.
From the previous proposition we have
$\Aut(\WS(\gt{u}),\cas)=\Sy(\gt{u})$. The group
$\Aut(\WS(\gt{u}),\wcgU)$ is just the subgroup of
$\Aut(\WS(\gt{u}))$ that fixes $\cej$. Consider the automorphism of
$\WS(\gt{u})$ obtained by setting $\ogi^{1/2}=\exp(\pi J(0)/2)$. It
is the automorphism of $\WS(\gt{u})$ induced by the orthogonal
transformation of $\gt{u}$ which is multiplication by $\ii$ on
$\gt{a}$, and multiplication by $-\ii$ on $\gt{a}^*$. Clearly, any
element of $\Aut(\WS(\gt{u}),\wcgU)$ commutes with $\ogi^{1/2}$. On
the other hand, if $g\in\Sy(\gt{u})$ commutes with $\ogi^{1/2}$ then
$g$ preserves the decomposition $\gt{u}=\gt{a}\oplus\gt{a}^*$, and
with respect to the basis $\{a_1,\ldots,a_1^*,\ldots\}$ is
represented by a block matrix
\begin{gather}\label{blckdecomp}
     g\sim
          \left(
              \begin{array}{cc}
                A & 0 \\
                0 & B \\
              \end{array}
            \right)
\end{gather}
and we must be able to write $A=T_g^t$ and $B=T_g^{-1}$ for some
invertible $N\times N$ matrix $T_g$ if it is the case that $g$
represents an element of $\Sy(\gt{u})$. (Here $T_g^t$ denotes the
transpose of $T_g$.) Evidently $\Aut(\WS(\gt{u}),\wcgU)$ is the
centralizer of $\ogi^{1/2}$ in $\Sy(\gt{u})$, and this group is all
matrices of the form $(\ref{blckdecomp})$, with $A$ invertible and
$B$ the inverse transpose of $A$. This is a copy of the group
$\GL(\gt{a})$.
\end{proof}

We record the results of this section in
\begin{thm}\label{thm:GLn:main}
For $N$ a positive integer, the quadruple $(\wu,Y,\vac,\wcgU)$ is a
self-dual $\U(1)$--VOA of rank $-N$. The full automorphism group of
$(\wu,\wcgU)$ is $\GL(\gt{a})$.
\end{thm}


\section{The Rudvalis group}\label{sec:Ru}

In this section we realize the sporadic simple group of Rudvalis as
symmetries of an enhanced VOA $(\wru,\wcgRu)$. Our plan for
realizing the Rudvalis group as symmetry of an enhanced VOA is to
consider first the enhanced VOA for ${\GL}_N(\CC)$ for suitable $N$,
constructed in \S\ref{sec:GLn}, and then to find a single extra
conformal generator with which to refine the conformal structure,
just to the point that the Rudvalis group becomes visible.

We give two constructions of this extra conformal generator. The
first construction arises directly from the geometry of the
Conway--Wales lattice \cite{ConRu}, and is given in
\S\ref{sec:Ru:geom}. We include also in this section an explicit
description of this lattice. The second construction, which is
perhaps more convenient for computations, is a description in terms
of a particular maximal subgroup of a double cover of the Rudvalis
group, and is given in \S\ref{sec:Ru:mnml}. Just as in the companion
article \cite{DunVARuI}, \S\ref{sec:Ru:geom} and \S\ref{sec:Ru:mnml}
are independent, and the reader my safely skip one in favor of the
other. The approach of \S\ref{sec:Ru:geom} is extremely simple.

The enhanced VOA structure for $(\wru,\wcgRu)$ is described in
\S\ref{sec:Ru:const}. We show in \S\ref{sec:Ru:symms} that the full
automorphism group of $(\wru,\wcgRu)$ is a certain four fold cover
of the Rudvalis group.

\subsection{Geometric description}\label{sec:Ru:geom}

The {\em Conway--Wales lattice} is a certain Hermitian lattice
$\LL_{\Ru}$ of rank $28$ over $\ZZ[\ii]$. Regarded as a lattice of
rank $56$ over $\ZZ$, it is even and self-dual, and has no roots
(elements of square length $2$). The full automorphism group of
$\LL_{\Ru}$ is a four fold cover of the sporadic simple group $\Ru$.

We now recall Conway's description of $\LL_{\Ru}$ from \cite{ConRu}.
We use a complex quaternionic notation which allows to regard a
typical vector in $\LL_{\Ru}$ as a seven tuple $(q_0,\cdots,q_6)$
where each $q_i$ is an element of the {\em complex quaternions}
which is the $\CC$-algebra with generators $\qj,\qk,\ql$ satisfying
the usual quaternion relations
\begin{gather}
     \qj\qk\ql=-1,\quad\qj\qk=-\qk\qj,\quad
          \qk\ql=-\ql\qk,\quad\&c.
\end{gather}
and with $\ii$ as usual denoting a square root of $-1$ in $\CC$ (so
that $\ii$ commutes with $\qj$, $\qk$, and $\ql$).

We take $\{\ii,\qj,\qk,\ql\}$ as a basis for the complex quaternions
as a vector space over $\CC$. We will write a vector
$\lambda\in\LL_{\Ru}$ as a $4\times 7$ array of complex numbers,
with the rows displaying our component quaternions, and the columns
denoting the $\CC$-coefficients of a given quaternion with respect
to our chosen basis elements $\ii,\qj,\qk,\ql$. Actually, all the
entries will lie in $\hZZ[\ii]$.

Regarding elements of $\CC^{28}$ as seven-tuples of complex
quaternions allows to describe easily the action of a certain group
of the shape $Q_8\times 2^3\!:\!7$ on $\CC^{28}$. The factor $Q_8$
is the quaternion group of order $8$ consisting of componentwise
right multiplications by the quaternions $\pm 1$, $\pm\qj$,
$\pm\qk$, $\pm\ql$. The $2^3$ of the right hand factor consists of
sign changes on coordinates $q_n$, $q_{n+3}$, $q_{n+5}$, $q_{n+6}$,
for $n\in\ZZ/7$ (with indices read modulo $7$), and the $7$ of the
right hand factor is the group of cyclic permutations of quaternion
components $(q_0,q_1,\ldots,q_6)\mapsto
(q_n,q_{1+n},\ldots,q_{6+n})$.

We now define $\LL_{\Ru}$ to be the $\ZZ[\ii]$-module generated by
the vectors of Tables~\ref{tab:Ru:geom:latvecs1} and
\ref{tab:Ru:geom:latvecs2}, and their images under the group
$Q_8\times 2^3\!:\!7$ just described. (In the tables we write
$\bar{a}$ as a shorthand for $-a$. The content of
Tables~\ref{tab:Ru:geom:latvecs1} and \ref{tab:Ru:geom:latvecs2} is
just Table~1 of \cite{ConRu}.)
\begin{table}
   \centering
   \caption{Some vectors in $\LL_{\Ru}$}
   \label{tab:Ru:geom:latvecs1}
\begin{small}
\begin{gather*}
\begin{array}{|rrrr|}
     2 & 0 & 0 & 0\\
     0 & 0 & 0 & 0\\
     0 & 0 & 0 & 0\\
     0 & 0 & 0 & 0\\
     0 & 0 & 0 & 0\\
     0 & 0 & 0 & 0\\
     0 & 0 & 0 & 0\\
\end{array}\\
\frac{1}{2}
\begin{array}{|rrrr|}
     \bar{\ii} & 1 & 1 & 1\\
     0 & 0 & 0 & 0\\
     0 & 0 & 0 & 0\\
     1 & 1 & \ii & \bar{1}\\
     0 & 0 & 0 & 0\\
     1 & \ii & \bar{1} & 1\\
     1 & \bar{1} & 1 & \ii\\
\end{array}\quad
\frac{1}{2}
\begin{array}{|rrrr|}
     \bar{\ii} & 1 & \bar{1} & \bar{1}\\
     0 & 0 & 0 & 0\\
     0 & 0 & 0 & 0\\
     1 & \bar{1} & \ii & 1\\
     0 & 0 & 0 & 0\\
     1 & \bar{\ii} & 1 & 1\\
     1 & \bar{1} & \bar{1} & \bar{\ii}\\
\end{array}\quad
\frac{1}{2}
\begin{array}{|rrrr|}
     \bar{\ii} & \bar{1} & 1 & \bar{1}\\
     0 & 0 & 0 & 0\\
     0 & 0 & 0 & 0\\
     1 & 1 & \bar{\ii} & 1\\
     0 & 0 & 0 & 0\\
     1 & \bar{\ii} & \bar{1} & \bar{1}\\
     1 & 1 & \bar{1} & \ii\\
\end{array}\quad
\frac{1}{2}
\begin{array}{|rrrr|}
     \bar{\ii} & \bar{1} & \bar{1} & 1\\
     0 & 0 & 0 & 0\\
     0 & 0 & 0 & 0\\
     1 & \bar{1} & \bar{\ii} & \bar{1}\\
     0 & 0 & 0 & 0\\
     1 & \ii & 1 & \bar{1}\\
     1 & 1 & 1 & \bar{\ii}\\
\end{array}
\end{gather*}
\end{small}
\end{table}

\begin{table}
   \centering
   \caption{Some more vectors in $\LL_{\Ru}$}
   \label{tab:Ru:geom:latvecs2}
\begin{small}
\begin{gather*}
\frac{1}{2}
\begin{array}{|rrrr|}
     1 & 1 & 1 & 1\\
     1 & 0 & 0 & \bar{1}\\
     1 & \bar{1} & 0 & 0\\
     0 & 0 & \bar{1} & 1\\
     1 & 0 & \bar{1} & 0\\
     0 & \bar{1} & 1 & 0\\
     0 & 1 & 0 & \bar{1}\\
\end{array}\quad
\frac{1}{2}
\begin{array}{|rrrr|}
     1 & 1 & 1 & 1\\
     1 & 0 & 0 & \bar{\ii}\\
     1 & \bar{\ii} & 0 & 0\\
     0 & 0 & \bar{\ii} & \bar{\ii}\\
     1 & 0 & \bar{\ii} & 0\\
     0 & \bar{\ii} & \bar{\ii} & 0\\
     0 & \bar{\ii} & 0 & \bar{\ii}\\
\end{array}\quad
\frac{1}{2}
\begin{array}{|rrrr|}
     1 & 1 & 1 & 1\\
     0 & \bar{1} & \ii & 0\\
     0 & 0 & \bar{1} & \ii\\
     \bar{1} & \ii & 0 & 0\\
     0 & \ii & 0 & \bar{1}\\
     \bar{1} & 0 & 0 & \ii\\
     \bar{1} & 0 & \ii & 0\\
\end{array}\quad
\frac{1}{2}
\begin{array}{|rrrr|}
     1 & 1 & 1 & 1\\
     0 & \ii & \bar{1} & 0\\
     0 & 0 & \ii & \bar{1}\\
     \ii & \bar{1} & 0 & 0\\
     0 & \bar{1} & 0 & \ii\\
     \ii & 0 & 0 & \bar{1}\\
     \ii & 0 & \bar{1} & 0\\
\end{array}\\
\frac{1}{2}
\begin{array}{|rrrr|}
     1 & 1 & \bar{1} & \bar{1}\\
     0 & 1 & \ii & 0\\
     0 & 0 & \ii & 1\\
     \bar{1} & \bar{\ii} & 0 & 0\\
     \bar{\ii} & 0 & \bar{\ii} & 0\\
     0 & \bar{\ii} & \bar{\ii} & 0\\
     \bar{\ii} & 0 & \bar{1} & 0\\
\end{array}\quad
\frac{1}{2}
\begin{array}{|rrrr|}
     1 & 1 & \bar{1} & \bar{1}\\
     0 & \bar{\ii} & \bar{1} & 0\\
     0 & 0 & \bar{1} & \bar{\ii}\\
     \ii & 1 & 0 & 0\\
     \bar{1} & 0 & 1 & 0\\
     0 & 1 & \bar{1} & 0\\
     \bar{1} & 0 & \bar{\ii} & 0\\
\end{array}\quad
\frac{1}{2}
\begin{array}{|rrrr|}
     1 & 1 & \bar{1} & \bar{1}\\
     1 & 0 & 0 & 1\\
     \bar{\ii} & \ii & 0 & 0\\
     0 & 0 & 1 & 1\\
     0 & \bar{1} & 0 & \ii\\
     \ii & 0 & 0 & \bar{1}\\
     0 & \ii & 0 & \bar{\ii}\\
\end{array}\quad
\frac{1}{2}
\begin{array}{|rrrr|}
     1 & 1 & \bar{1} & \bar{1}\\
     \bar{\ii} & 0 & 0 & \ii\\
     1 & 1 & 0 & 0\\
     0 & 0 & \ii & \bar{\ii}\\
     0 & \bar{\ii} & 0 & 1\\
     1 & 0 & 0 & \bar{\ii}\\
     0 & 1 & 0 & 1\\
\end{array}\\
\frac{1}{2}
\begin{array}{|rrrr|}
     1 & \bar{1} & 1 & \bar{1}\\
     \bar{\ii} & 0 & 0 & \bar{\ii}\\
     0 & 0 & 1 & \ii\\
     0 & 0 & \bar{\ii} & \bar{\ii}\\
     0 & 1 & 0 & \ii\\
     \bar{\ii} & 0 & 0 & \bar{1}\\
     \bar{1} & 0 & \bar{\ii} & 0\\
\end{array}\quad
\frac{1}{2}
\begin{array}{|rrrr|}
     1 & \bar{1} & 1 & \bar{1}\\
     \bar{1} & 0 & 0 & 1\\
     0 & 0 & \bar{\ii} & \bar{1}\\
     0 & 0 & 1 & \bar{1}\\
     0 & \bar{\ii} & 0 & \bar{1}\\
     \bar{1} & 0 & 0 & \bar{\ii}\\
     {\ii} & 0 & 1 & 0\\
\end{array}\quad
\frac{1}{2}
\begin{array}{|rrrr|}
     1 & \bar{1} & 1 & \bar{1}\\
     0 & \ii & \bar{1} & 0\\
     1 & 1 & 0 & 0\\
     \ii & \bar{1} & 0 & 0\\
     \bar{\ii} & 0 & \ii & 0\\
     0 & \bar{\ii} & \ii & 0\\
     0 & 1 & 0 & 1\\
\end{array}\quad
\frac{1}{2}
\begin{array}{|rrrr|}
     1 & \bar{1} & 1 & \bar{1}\\
     0 & 1 & \bar{\ii} & 0\\
     \bar{\ii} & \ii & 0 & 0\\
     1 & \bar{\ii} & 0 & 0\\
     1 & 0 & 1 & 0\\
     0 & 1 & 1 & 0\\
     0 & \bar{\ii} & 0 & \ii\\
\end{array}\\
\frac{1}{2}
\begin{array}{|rrrr|}
     1 & \bar{1} & \bar{1} & 1\\
     0 & \ii & 1 & 0\\
     \bar{\ii} & \bar{\ii} & 0 & 0\\
     \bar{\ii} & \bar{1} & 0 & 0\\
     0 & \ii & 0 & 1\\
     \bar{1} & 0 & 0 & \bar{\ii}\\
     0 & \bar{\ii} & 0 & \bar{\ii}\\
\end{array}\quad
\frac{1}{2}
\begin{array}{|rrrr|}
     1 & \bar{1} & \bar{1} & 1\\
     0 & \bar{1} & \bar{\ii} & 0\\
     \bar{1} & 1 & 0 & 0\\
     \bar{1} & \bar{\ii} & 0 & 0\\
     0 & \bar{1} & 0 & \bar{\ii}\\
     \ii & 0 & 0 & 1\\
     0 & \bar{1} & 0 & 1\\
\end{array}\quad
\frac{1}{2}
\begin{array}{|rrrr|}
     1 & \bar{1} & \bar{1} & 1\\
     \bar{\ii} & 0 & 0 & \ii\\
     0 & 0 & \ii & \bar{1}\\
     0 & 0 & \bar{\ii} & \ii\\
     1 & 0 & 1 & 0\\
     0 & 1 & 1 & 0\\
     \ii & 0 & \bar{1} & 0\\
\end{array}\quad
\frac{1}{2}
\begin{array}{|rrrr|}
     1 & \bar{1} & \bar{1} & 1\\
     1 & 0 & 0 & 1\\
     0 & 0 & 1 & \bar{\ii}\\
     0 & 0 & 1 & 1\\
     \bar{\ii} & 0 & \ii & 0\\
     0 & \ii & \bar{\ii} & 0\\
     1 & 0 & \bar{\ii} & 0\\
\end{array}
\end{gather*}
\end{small}
\end{table}

Now set $\gt{r}=\CC\otimes_{\ZZ[\ii]} \LL_{\Ru}$ and
$\gt{s}=\gt{r}\oplus\gt{r}^*$, and consider the element $\ossp\in
\bigvee^4(\gt{r})\oplus\bigvee^4(\gt{r}^*)$ described as
\begin{gather}
     \ossp=\sum_{\lambda\in(\LL_{\Ru})_2}\lambda^4
          +(\lambda^*)^4
\end{gather}
where we write $(\LL_{\Ru})_2$ for the type $2$ (i.e. square norm
$4$) vectors in $\LL_{\Ru}$.
\begin{prop}
The element $\ossp$ is non-zero.
\end{prop}
\begin{proof}
The vector $\sum_{(\LL_{\Ru})_2}(\lambda^*)^4$ may be viewed as a
polynomial function on $\gt{r}$. Consider the value of this function
at the first vector of Table \ref{tab:Ru:geom:latvecs1} say, so that
we are just computing the sum of the fourth powers of the top left
entries of each minimal length vector in $\LL_{\Ru}$. Since all
entries lie in $\hZZ\cup \hZZ\ii$, each term in the sum is
non-negative. At least one term is positive. We conclude that
$\ossp$ is non-zero.
\end{proof}

Note that the space $\bigvee^n(\gt{r})$ embeds naturally in
$\WS(\gt{s})$ under the map
\begin{gather}
     a_{i_1}\vee\cdots\vee a_{i_n}\mapsto
     a_{i_1}(-\tfrac{1}{2})\cdots
          a_{i_n}(-\tfrac{1}{2})\vac
\end{gather}
and similarly for $\bigvee^n(\gt{r}^*)$. In \S\ref{sec:Ru:const} we
will identify $\ossp$ with its image in $\WS(\gt{s})_2$ under this
embedding.

\subsection{Monomial description}\label{sec:Ru:mnml}

In this section we consider the invariants of the monomial group $M$
in the spaces $\bigvee^4(\gt{r})$ and $\bigvee^4(\gt{r}^*)$, which
both embed naturally in $\WS(\gt{s})_2$. We will assume that the
reader is acquainted with the notations and terminology of \S5.2 of
the preceding article \cite{DunVARuI}.

\medskip

Consider first the space $\bigvee^4(\gt{r})$. A basis for this space
is given by the homogeneous monomials of degree $4$ in the symbols
$\{a_i\}_{i\in\hfset}$, and the embedding in $\WS(\gt{s})_2$ is
given by $f(a_i)\mapsto f(a_i(-\tfrac{1}{2}))\vac$, for $f$ a
monomial (homogeneous of degree $4$).

Recall the code $\mc{D}$ generated by the dozens of $\hfset$.
Suppose that $t\in\bigvee^4(\gt{r})$ is invariant for the action of
$M$, and let us consider the expansion of $t$ with respect to the
basis of monomials of degree $4$. Then it is easy to see that the
coefficient of a monomial of the form $a_i^3a_j$ must be zero for
all $\{i,j\}\subset\hfset$, since for each such pair in $\hfset$
there is an element of $A\subset M$ that fixes $a_i$ say, and
negates $a_j$; in other words, the dual code to $\mc{D}$ has no
words of weight $2$. Similarly, the coefficient of any monomial of
the from $a_{i}^2a_{jk}$ must also vanish. The coefficient of a
monomial of the form $a_I$ with $I\subset\hfset$ of cardinality $4$
must be zero unless $I$ lies in the dual to $\mc{D}$, and the words
of weight four in $\mc{D}^*$ are precisely the quartets of $\hfset$.
The ringed quartets $\mc{R}$ and the non-ringed quartets $\mc{N}$
constitute the two orbits of $\overline{M}$ on weight four words in
$\mc{D}^*$.

We see that there are four kinds of monomials that can appear with
non-trivial coefficient in an invariant for the action of $M$ on
$\bigvee^4(\gt{r})$. Namely, those of the form $a_i^4$ for
$i\in\hfset$, or $a_C^2$ for $C\in\mc{C}$ a couple in $\hfset$, or
$a_N$ for $N\in\mc{N}$, or $a_R$ for $R\in\mc{R}$. Since
$\overline{M}$ is transitive on each of the sets $\hfset$, $\mc{C}$,
$\mc{R}$, and $\mc{N}$, it follows that there is at most a
four-space of invariants for the action of $M$ on
$\bigvee^4(\gt{r})$. Candidates for a spanning set for this space
are the vectors
\begin{gather}\label{eqn:Ru:monomialinvs:tis}
     t_1=\sum_{i\in\hfset}\varepsilon_ia_i^4,\quad
     t_2=\sum_{C\subset\mc{C}}\varepsilon_C a_C^2,\quad
     t_3=\sum_{R\subset\mc{R}}\varepsilon_Ra_R,\quad
     t_4=\sum_{N\subset\mc{N}}\varepsilon_Na_N
\end{gather}
where the functions $X\mapsto\varepsilon_X$ are uniquely determined
up to scalar multiplication, so long as they exist; that is, give
rise to a non-zero invariant. Recall that each element of $M$ may be
written as a coordinate permutation followed by a diagonal matrix
with entries in $\{\pm 1,\pm \ii\}$. It follows that we may set
$\varepsilon_i=1$ for all $i$ in $\hfset$, and that $t_1$ is then
invariant for the action of $M$. More generally, we have
\begin{prop}\label{prop:Ru:monomialinvs:cocycles}
Set $\varepsilon_i=1$ for each $i\in\hfset$, and $\varepsilon_R=1$
for each $R\in\mc{R}$. Set $\varepsilon_C=1$ for each couple $C$
associated to a dozen containing $\infty$, and $\varepsilon_C=-1$
otherwise. Set $\varepsilon_N=-1$ for each non-ringed quartet $N$
that does not contain $\infty$ and has intersection at most one with
each block, and $\varepsilon_N=1$ otherwise. Then the vectors $t_1$,
$t_2$, $t_3$, and $t_4$ of (\ref{eqn:Ru:monomialinvs:tis}), are each
invariant for the action of $M$ on $\bigvee^4(\gt{r})$.
\end{prop}
Note that the functions $\varepsilon$ may all be chosen to take
values in $\{\pm 1\}$. We thus obtain four linearly independent
invariants $t_i^*$ for the action of $M$ on $\bigvee^4(\gt{r}^*)$ by
replacing $a_i$ with $a_i^*$ in the expressions of
(\ref{eqn:Ru:monomialinvs:tis}).

We identify the $t_i$ and $t_i^*$ with their images in
$\WS(\gt{s})_2$, so that $t_1=\sum a_i(-\tfrac{1}{2})^4\vac$, and
$t_1^*=\sum a_i^*(-\tfrac{1}{2})^4\vac$, \&c., and we define a
vector $\ossp\in\WS(\gt{s})_2$ by setting
\begin{gather}
\begin{split}
     c_0\ossp=&\;13t_1
          +(18+12\ii)t_2
          +(-192-24\ii)t_3
          +(-48+72\ii)t_4\\
          &+13t_1^*
          +(18-12\ii)t_2^*
          +(-192+24\ii)t_3^*
          +(-48-72\ii)t_4^*
\end{split}
\end{gather}
where $c_0=2^7.3.5.13.29$.

\subsection{Construction}\label{sec:Ru:const}

Recall the VOA $\WS(\gt{u})$ of \S\ref{sec:weylalgs:VOAs}, and the
vector space $\gt{s}$ of \S\ref{sec:Ru:geom}. We now define $\wru$
to be the enhanced VOA with underlying VOA structure given by
$\wru=\WS(\gt{s})$, and with enhanced conformal generators given by
$\wcgRu=\{\cas,\cej,\ossp\}$, where $\delta$ is as in
\S\ref{sec:Ru:geom} or \S\ref{sec:Ru:mnml} (the descriptions are
equivalent), and $\cas$ and $\cej$ are as in \S\ref{sec:GLn}. We
have
\begin{thm}\label{thm:Ru:const:eVOA}
The quadruple $(\wru,Y,\vac,\wcgRu)$ is a self-dual enhanced
$\U(1)$--VOA of rank $-28$.
\end{thm}

\subsection{Symmetries}\label{sec:Ru:symms}

We now consider the automorphism group of the enhanced VOA $\wru$
with enhanced conformal structure given by $\wcgRu$. We opt to use
the construction of \S\ref{sec:Ru:geom} for the purposes of
describing $\Aut(\wru)$.

Let $F$ be the subgroup of ${\GL}(\gt{r})$ that fixes $\ossp$. Then
$F$ is the full automorphism group of the enhanced VOA structure on
$\wru$. The element $\ossp$ is evidently invariant under
$\Aut(\LL_{\Ru})$, and thus we have an embedding
$4.\!\Ru\hookrightarrow F$.

\begin{prop}\label{prop:Ru:symms:FhasR}
The group $F$ contains a group isomorphic to $4.\!\Ru$.
\end{prop}

We write $R$ for the copy of $4.\!\Ru$ in ${\GL}(\gt{r})$. The proof
of the following proposition is directly parallel to the proof of
Proposition 5.10 of Part I \cite{DunVARuI}. The main ingredient is
the fact that the only elements of the Lie algebra of
${\GL}(\gt{r})$ that annihilate $\ossp$ are those that exponentiate
to scalar matrices. The nature of $\ossp$ then forces these scalars
to be fourth roots of unity.
\begin{prop}\label{prop:Ru:symms:Fisfinite}
The group $F$ is finite, and is contained in $\SL(\gt{r})$.
\end{prop}
We now apply Proposition 5.12 of Part I \cite{DunVARuI} to conclude
that the group $F$ is precisely the group $4.\!\Ru$. Actually, the
Schur multiplier of the Rudvalis group has order $2$ \cite{ATLAS};
the center of $F$ is evidently generated by scalar multiplications
by $\ii$, and the simple group $\Ru$ has no non-trivial
representation of degree $28$, so the only possibility is that the
group $4.\!\Ru$ is in fact a central product of this cyclic group
$\lab\,\ii \Id\rab$ with the perfect double cover of $\Ru$.

We record this and the other results of this section in the
following
\begin{thm}\label{thm:Ru:main}
The quadruple $(\wru,Y,\vac,\wcgRu)$ is a self-dual enhanced
$\U(1)$--VOA of rank $-28$. The full automorphism group of
$(\wru,\wcgRu)$ is a central product of a cyclic group of order four
with the perfect double cover of the sporadic simple group of
Rudvalis.
\end{thm}


\section{McKay--Thompson series}\label{sec:series}

In this section we consider the McKay--Thompson series arising from
the enhanced VOA constructed in \S\ref{sec:Ru}. These series furnish
an analogue of Monstrous Moonshine for the sporadic group of
Rudvalis, and we will derive explicit expressions for them here.

The main tool for expressing the series explicitly is the notion of
weak Frame shape, discussed in the companion article
\cite[\S6.3]{DunVARuI}. We refer to there for the definition of weak
Frame shape.

We recall briefly the notion of two variable McKay--Thompson series
associated to the action of an element of the automorphism group of
an enhanced $\U(1)$--VOA in \S\ref{sec:series:2var}, and we then
present explicit expressions for all the McKay--Thompson series
arising from the action of $2.\!\Ru$ on $\wru$.

As mentioned in \S\ref{sec:intro}, we obtain a genus zero property
for the Rudvalis group by considering together the series arising
here and those of the preceding article \cite{DunVARuI}. We define
the functions $\tilde{F}^A_g(\tau)$ and
$\tilde{F}_{g}^{\forall}(\tau)$ for $g\in \Ru$, and make the main
observation regarding these functions in \S\ref{sec:series:MbeyM}.

\subsection{Two variable McKay--Thompson
series}\label{sec:series:2var}

Let $(U,Y,\vac,\{\cas,\cej,\ldots\})$ be an enhanced $\U(1)$--VOA of
rank $c$ (see \S\ref{sec:eVOAs}). Then for $g$ an automorphism of
$(U,\{\cas,\cej,\ldots\})$, the {\em two variable McKay--Thompson
series associated to the action of $g$ on $U$} is the series in
variables $p$ and $q$ given by
\begin{gather}
     \tr_{U}gp^{J(0)}q^{L(0)-c/24}
     =\sum_{m,n}(\tr_{U_{n}^m}g)p^mq^{n-c/24}
\end{gather}
where $U_n^m$ denotes the subspace of $U$ of degree $n$ consisting
of vectors of charge $m$ (see \S\ref{sec:eVOAs}). In the limit as
$p\to 1$ we recover the (ordinary) McKay--Thompson series associated
to the action of $g$ on $U$. In the case that $g$ is the identity we
obtain what we call the {\em two variable character of $U$}.

We would like to compute the two variable McKay--Thompson series
arising from our main example $\wru$. In this case we are in the
situation of \S\ref{sec:GLn}, so that the relevant automorphisms lie
in $\GL(\gt{a})$ where $\gt{a}$ is a Hermitian vector space, and the
VOA underlying the relevant enhanced VOA is of the form
$\WS(\gt{u})$ for $\gt{u}=\gt{a}\oplus \gt{a}^*$.

We suppose that $g$ lies in $\GL(\gt{a})$, and has a weak Frame
shape $\prod_j(k_j)_{a_j}^{m_j}$ say, for its action on $\gt{a}$.
Recall the Jacobi theta function $\vartheta(z|\tau)$ from
\S\ref{sec:intro:notation}. Then we set
\begin{align}\label{eqn:series:2var|phipsidef}
     \phi_{g}(z|\tau)
          &=\prod_j
          \frac{\vartheta(k_jz+1/2+a_j|k_j\tau)^{m_j}}
               {\eta(k_j\tau)^{m_j}}\\
     \psi_{g}(z|\tau)
          &=\prod_j
          p^{k_jm_j/2}q^{k_jm_j/8}
          \frac{\vartheta(k_jz+1/2+a_j+k_j\tau/2|k_j\tau)^{m_j}}
               {\eta(k_j\tau)^{m_j}}
\end{align}
(just as in \cite[\S6]{DunVARuI}) and we then have
\begin{thm}\label{thm:series:2var|chars}
Let $g\in\GL(\gt{a})$. Then the McKay--Thompson series associated to
the actions of $g$ and $\ogi g$ on $\WS(\gt{u})$ and
$\WS(\gt{u})_{\ogi}$ admit the following expressions.
\begin{align}\label{eqn:series:2var|chars}
    \mathsf{tr}|_{\WS(\gt{u})}gp^{J(0)}q^{L(0)-c/24}=
          &\,\frac{1}{\phi_{{g}}(z|\tau)}\\
    \mathsf{tr}|_{\WS(\gt{u})}\ogi gp^{J(0)}q^{L(0)-c/24}=
          &\,\frac{1}{\phi_{{-g}}(z|\tau)}\\
    \mathsf{tr}|_{\WS(\gt{u})_{\ogi}}gp^{J(0)}q^{L(0)-c/24}=
          &\,\frac{1}{\psi_{{g}}(z|\tau)}\\
    \mathsf{tr}|_{\WS(\gt{u})_{\ogi}}\ogi gp^{J(0)}q^{L(0)-c/24}=
          &\,\frac{1}{\psi_{{-g}}(z|\tau)}
\end{align}
\end{thm}
The verification of Theorem \ref{thm:series:2var|chars} is very
similar to that of the corresponding Theorem 6.1 in \cite{DunVARuI}.

\medskip

Recall that the vector space underlying the enhanced vertex algebra
$\aru$ constructed in Part I \cite{DunVARuI} may be expressed as a
direct sum $\aru=A(\gt{s})^0\oplus A(\gt{s})^0_{\ogi}$. (Where the
$\gt{s}$ there and in this article coincide as Hermitian vector
spaces.) The action of $\Ru$ preserves this decomposition, and more
than this, $\Ru$ acts projectively on the companion spaces
$A(\gt{s})^1$ and $A(\gt{s})^1_{\ogi}$. In fact, we may regard the
object $A(\gt{s})=A(\gt{s})^0\oplus A(\gt{s})^1$ as an enhanced
$\U(1)$--VOA with a projective action by $\Ru$, and
$A(\gt{s})_{\ogi}$ is then its canonically twisted module, again
with a projective action by $\Ru$. We now set $\arutd=A(\gt{s})$ and
$(\arutd)_{\ogi}=A(\gt{s})_{\ogi}$. We regard $\arutd$ as a twisting
of $\aru$, and $(\arutd)_{\ogi}$ is then its canonically twisted
module. This twisting breaks much of the enhanced conformal
structure on $\aru$, but the $\U(1)$--VOA structure at least
survives.

For convenience later we now record the analogue of Theorem
\ref{thm:series:2var|chars} with $A(\gt{u})$ and $A(\gt{u})_{\ogi}$
in place of $\WS(\gt{u})$ and $\WS(\gt{u})_{\ogi}$.
\begin{thm}\label{thm:series:2var|CMchars}
Let $g\in\GL(\gt{a})$. Then the McKay--Thompson series associated to
the actions of $g$ and $\ogi g$ on $A(\gt{u})$ and
$A(\gt{u})_{\ogi}$ admit the following expressions.
\begin{align}\label{eqn:series:2var|CMchars}
    \mathsf{tr}|_{A(\gt{u})}gp^{J(0)}q^{L(0)-c/24}=
          &\,{\phi_{{-g}}(z|\tau)}\\
    \mathsf{tr}|_{A(\gt{u})}\ogi gp^{J(0)}q^{L(0)-c/24}=
          &\,{\phi_{{g}}(z|\tau)}\\
    \mathsf{tr}|_{A(\gt{u})_{\ogi}}gp^{J(0)}q^{L(0)-c/24}=
          &\,{\psi_{{-g}}(z|\tau)}\\
    \mathsf{tr}|_{A(\gt{u})_{\ogi}}\ogi gp^{J(0)}q^{L(0)-c/24}=
          &\,{\psi_{{g}}(z|\tau)}
\end{align}
\end{thm}

\subsection{Moonshine beyond the Monster}\label{sec:series:MbeyM}

Using Theorem \ref{thm:series:2var|chars}, we may now compute the
two variable McKay--Thompson series associated to the action of any
element $g\in 4.\!\Ru$ on $\wru$ as soon as the corresponding weak
Frame shape is known. We reproduce in Table~\ref{tab:series:Ru|frms}
the Frame shapes and weak Frame shapes for the Rudvalis group, given
also in \cite[\S6]{DunVARuI}.
\begin{table}
  \centering
  \caption{Frame Shapes for $2.\!\Ru$}
  \label{tab:series:Ru|frms}
\renewcommand{\arraystretch}{1.1}
\begin{tabular}{c|cc||c|cc}
  Class & $\SO_{56}$ & ${\SU}_{28}$ & Class & $\SO_{56}$ & ${\SU}_{28}$   \\  \hline
    1A &  $1^{56}$ &  $1^{28}$ &12B & $4^112^5/2^16^1$ & $1_{1/4}^13_{3/4}^112^2$   \\
    2A &  $1^82^{24}$ &    $1^42^{12}$ &13A & $1^413^4$ & $1^213^2$ \\
    2B &  $4^{28}/2^{28}$ &    $4^{14}/2^{14}$ &14A & $28^4/14^4$ & $28^2/14^2$ \\
    3A &  $1^23^{18}$ &  $1^13^9$ &14B & $28^4/14^4$ & $28^2/14^2$  \\
    4A &  $1^84^{12}$ &  $1^44^6$ &14C & $28^4/14^4$ & $28^2/14^2$  \\
    4B &  $4^{16}/2^4$ & $1_{1/4}^44^{6}$ &15A & $1^215^4/3^2$ & $1^115^2/3^1$  \\
    4C &  $4^{16}/2^4$ & $4^8/2^2$ &16A & $1^216^4/2^18^1$ & $1^11_{1/4}^12_{3/4}^116^2/8^1$\\
    4D &  $2^44^{12}$ &  $2^24^6$ &16B & $1^216^4/2^18^1$ & $1^11_{3/4}^12_{1/4}^116^2/8^1$ \\
    5A &  $1^65^{10}$ &  $1^35^5$ &20A & $1^24^210^220^2/2^25^2$ & $1^14^110^120^1/2^15^1$ \\
    5B &  $5^{12}/1^4$ & $5^{6}/1^2$ &20B & $2^14^120^3/10^1$ & $1_{1/4}^12^15_{3/4}^120^1$ \\
    6A &  $1^23^26^8$ & $1^13^16^4$ &20C & $2^14^120^3/10^1$ & $1_{1/4}^12^15_{3/4}^120^1$\\
    7A &  $7^8$ &   $7^4$ &24A & $4^112^124^2/2^16^1$ & $1_{1/4}^13_{1/4}^124^1$ \\
    8A &  $4^48^6/2^4$ & $1_{1/4}^24^18^3/2^1$ &24B & $4^112^124^2/2^16^1$ & $1_{1/4}^13_{1/4}^124^1$ \\
    8B &  $8^8/4^2$ & $1_{1/4}^22^18^4/4^2$ &26A & $4^252^2/2^226^2$ & $4^152^1/2^126^1$\\
    8C &  $4^28^6$ & $4^18^3$ &26B & $4^252^2/2^226^2$ & $4^152^1/2^126^1$ \\
    10A & $2^45^210^4/1^2$ & $2^25^110^2/1^1$ &26C & $4^252^2/2^226^2$ & $4^152^1/2^126^1$\\
    10B & $2^220^6/4^210^6$ & $2^120^3/4^110^3$ &29A & $29^2/1^2$ & $29^1/1^1$ \\
    12A & $1^23^212^4$ & $1^13^112^2$ &29B & $29^2/1^2$ & $29^1/1^1$   \\
\end{tabular}
\end{table}

In Table~\ref{tab:osmt1A} we tabulate the terms of lowest charge and
degree in the character of $\wru$. The column headed $m$ is the
coefficient of $p^m$ (as a series in $q$), and the row headed $n$ is
the coefficient of $q^{n}$ (as a series in $p$). The coefficients of
$p^{-m}$ and $p^m$ coincide, and the coefficient of $p^mq^n$
vanishes unless $2m$ and $n$ are of the same parity.
\begin{table}
   \centering
     \caption{Character of $\wru$}
   \label{tab:osmt1A}
   \begin{tabular}{c|rrrrr}
     & 0 & 2 & 4 & 6 & 8 \\ \hline
     0& $1$&  &  &  &  \\
     1& $784$& $406$&  &  &  \\
     2& $166404$& $114464$& $31465$&  &  \\
     3& $17122560$& $13207964$& $5752208$& $1107568$&  \\
     4& $1083938457$& $889479360$& $480258212$& $156267776$& $23535820$\\
     5& $48023166576$& $40957573860$& $25100492032$& $10592918798$& $2786641088$\\
     6& $1612815529556$& $1412141567872$& $941094824285$& $466052603296$& $163110431792$\\
     7& $43324776509184$& $38661856383304$& $27348457208240$& $15108948493100$& $6331332324352$\\
   \end{tabular}\\
   \begin{tabular}{c|rrrr}
     & 1 & 3 & 5 & 7 \\ \hline
     1/2& $28$&  &  &  \\
     3/2& $11396$& $4060$&  &  \\
     5/2& $1681708$& $892388$& $201376$&  \\
     7/2& $135613016$& $85845452$& $31892924$& $5379616$\\
     9/2& $7178616060$& $5024198872$& $2375544164$& $690014864$\\
     11/2& $276842870136$& $206737356700$& $113074057180$& $43211564144$\\
     13/2& $8314619588660$& $6496652289400$& $3919853988680$& $1777585597900$\\
     15/2& $203616821215276$& $164493601241508$& $106534952459628$& $54381423797540$\\
   \end{tabular}
\end{table}
The entries of Table~\ref{tab:osmt1A} extend the coincidences
(\ref{eqn:intro_RuIIObs}) noted in \S\ref{sec:intro}, and in
addition, we now see that small degree characters of the covering
group $2.\!\Ru$ participate in similar identities also.
\begin{gather}\label{eqn:series:MbeyM|MTobsod}
     \begin{split}
     28&=28\\
     4060&=28+4032\\
     11396&=(2)28+4032+7308\\
     201376&=28+(2)4032+7308+87696+98280
     \end{split}
\end{gather}


We now focus attention on the series arising from the projective
action of the Rudvalis group on $(\arutd)_{\ogi}$ and
$(\wru)_{\ogi}$. For an arbitrary rational vertex operator
superalgebra (satisfying some technical conditions) the characters
of ordinary modules are not usually invariant for the action of the
modular group $\PSL_2(\ZZ)$, but rather for the index $3$ subgroup
generated by $\tau\mapsto \tau+2$ and $\tau\mapsto-1/\tau$. However,
considering the super characters of canonically twisted modules one
does expect to recover a representation of $\PSL_2(\ZZ)$; indeed,
this has been shown in \cite{DonZhaMdltyOrbVOSA}. By a similar
token, we expect the best modular properties to be enjoyed by the
series associated to the Rudvalis group via its action on the
twisted modules $(\arutd)_{\ogi}$ and $(\wru)_{\ogi}$.

The subspace of $(\arutd)_{\ogi}$ of lowest degree is a copy of the
exterior algebra $\bigwedge(\gt{r})$. We record the traces of
elements of $2.\!\Ru$ on this $2.\!\Ru$-module in
Table~\ref{tab:series:Ru|wedgechar}. A class determined by an
element $\hat{g}$ such that $\hat{g}$ and $-\hat{g}$ both have
vanishing trace on $\bigwedge(\gt{r})$ is omitted from the table. We
also record in Table~\ref{tab:series:Ru|wedgechar} the fusion of
conjugacy classes under the natural map $2.\!\Ru\to \Ru$, but only
for the classes in $\Ru$ that have a lift to $2.\!\Ru$ with
non-trivial trace on $\bigwedge(\gt{r})$. Our naming of conjugacy
classes follows the computer system \cite{GAP4}.

\comment{\begin{table}
  \centering
  \caption{Character for the action of $2.\!\Ru$ on $\bigwedge(\gt{r})$}
  \label{tab:series:Ru|wedgechar}
\begin{gather*}
\begin{array}[t]{c||cc|c|cc|cc|cc|cc|c|cc}
  \Ru & &{1A} & 2B & &{3A} && 5A && 5B && 7A & 10B && 13A  \\  \hline
  2.\!\Ru & 1A & 2A & 4A & 3A & 6A & 5A & 10A & 5B & 10B & 7A & 14A & 20A & 13A & 26A \\
  \hline
  \chi & 2^{28} & 0 & 2^{14} & 2^{10} &0 & 2^8 & 0 & 16 & 0 & 16 &
  0 & 4 & 16 & 0
\end{array}\\
\begin{array}[t]{c||c|c|c|cc|c|c|c|cc|cc}
  \Ru & 14A & 14B & 14C &&{15A} & 26A & 26B & 26C && 29A && 29B \\\hline
  2.\!\Ru & 28A & 28B & 28C & 15A & 30A & 52A & 52B & 52C & 29A & 58A & 29B & 58B \\ \hline
  \chi & 4& 4 & 4& 4&0 & 4&4&4& 1&29&1&29
\end{array}
\end{gather*}
\end{table}}

\begin{table}
  \centering
  \caption{Character for the action of $2.\!\Ru$ on $\bigwedge(\gt{r})$}
  \label{tab:series:Ru|wedgechar}
\renewcommand{\arraystretch}{1.1}
\begin{tabular}[t]{c||cc|c|cc|cc|cc|cc|c|cc}
  $\Ru$ & &{1A} & 2B & &{3A} && 5A && 5B && 7A & 10B && 13A  \\  \hline
  $2.\!\Ru$ & 1A & 2A & 4A & 3A & 6A & 5A & 10A & 5B & 10B & 7A & 14A & 20A & 13A & 26A \\
  \hline
  $\chi$ & $2^{28}$ & 0 & $2^{14}$ & $2^{10}$ &0 & $2^8$ & 0 & 16 & 0 & 16 &
  0 & 4 & 16 & 0\\
\end{tabular}
\begin{tabular}[t]{c||c|c|c|cc|c|c|c|cc|cc}
  $\Ru$ & 14A & 14B & 14C &&{15A} & 26A & 26B & 26C && 29A && 29B \\\hline
  $2.\!\Ru$ & 28A & 28B & 28C & 15A & 30A & 52A & 52B & 52C & 29A & 58A & 29B & 58B \\ \hline
  $\chi$ & 4& 4 & 4& 4&0 & 4&4&4& 1&29&1&29\\
\end{tabular}
\end{table}

We now define functions $\tilde{F}_g^A(\tau)$ and
$\tilde{F}_{g}^{\forall}(\tau)$ for $g\in \Ru$ by setting
\begin{gather}
     \tilde{F}_g^{A}(\tau)=
          \lim_{p\to 1}\tr_{(\arutd)_{\ogi}}
          -\hat{g}p^{J(0)}q^{L(0)-c/24}\label{eqn:series:Ftdsp}\\
     \tilde{F}_{g}^{\forall}(\tau)=
          \lim_{p\to 1}\tr_{(\wru)_{\ogi}}
          \hat{g}p^{J(0)}q^{L(0)-c/24}\label{eqn:series:Fos}
\end{gather}
where $\hat{g}$ is a preimage of $g$ in $2.\!\Ru$ of minimal order.
One may check by consulting the character table of $2.\!\Ru$ that
any two preimages of minimal order have the same Frame shape. We
define a vector space $\gt{F}_g$, for each $g\in \Ru$, by setting
\begin{gather}
     \gt{F}_g=\Span_{\CC}\left\{\tilde{F}_g^A(\tau),
          \tilde{F}_{g}^{\forall}(\tau)\right\}.
\end{gather}
In the case that $-\hat{g}$ has trace $0$ on $\bigwedge(\gt{r})$,
this amounts to setting $\tilde{F}_g^{A}(\tau)=0$ and
$\tilde{F}_{g}^{\forall}(\tau)=\infty$, and we invite the reader to
regard the corresponding space $\gt{F}_g$ as a copy of the trivial
two dimensional representation of $\PSL_2(\ZZ)$.

Care must be applied when taking the limit in
(\ref{eqn:series:Fos}). In the case of (\ref{eqn:series:Ftdsp}) each
coefficient in $q$ is a polynomial in $p$, and there is no problem
regarding the a priori formal variable $p$ as a variable on
$\CC^{\times}$. In the case of (\ref{eqn:series:Fos}) there is a
factor of the form
\begin{gather}\label{eqn:series:pprod}
     \prod_i\frac{1}{(1-\xi_ip)}
\end{gather}
where $\xi_i$ ranges over the eigenvalues of $\hat{g}$ (c.f. Theorem
\ref{thm:series:2var|chars}), which is formally a power series in
$p$, and this series converges only for $|p|<1$. In order to make
sense of the limit in (\ref{eqn:series:Fos}) we should analytically
continue each coefficient in $q$ as a function in $p$, so that the
series represented by (\ref{eqn:series:pprod}) may be replaced with
the product expression there --- in other words, the factor
(\ref{eqn:series:pprod}) should be regarded as the reciprocal of the
trace of $-\hat{g}$ on $\bigwedge(\gt{r})$.

From Table~\ref{tab:series:Ru|wedgechar} we see that $\gt{F}_g$ is
comprised of non-constant functions just in the case that $g$
belongs to one of the following $10$ conjugacy classes.
\begin{gather}
     2B\quad 10B\quad
     14A\quad 14B\quad 14C\quad
     26A\quad 26B\quad 26C\quad
     29A\quad 29B
\end{gather}
We will examine these cases in a some detail in the following
sections.

We refer the reader to \cite{ConNorMM}, \cite{ConMcKSebDiscGpsM} and
\cite{Cum_CngrSbGpsGenus01} for more information about discrete
subgroups of $\PSL_2(\RR)$ of genus zero. We adopt the notation of
\cite{ConNorMM} for discrete subgroups of $\PSL_2(\RR)$.

\subsubsection{$\gt{F}_{2B}$}

An element $g$ of class $2B$ in $\Ru$ lifts to a unique class of
order $4$ in the covering group $2.\!\Ru$, and an element of this
class has trace $2^{14}$ on $\bigwedge(\gt{r})$. We see from
Table~\ref{tab:series:Ru|frms} that the Frame shape of $\hat{g}$ is
$4^{28}/2^{28}$, so we have
\begin{gather}
     \tilde{F}_{2B}^A(\tau)=\frac{2^{14}}{1}
          \frac{\eta(4\tau)^{28}}{\eta(2\tau)^{28}},\quad
     \tilde{F}_{2B}^{\forall}(\tau)=\frac{1}{2^{14}}
          \frac{\eta(2\tau)^{28}}{\eta(4\tau)^{28}}.
\end{gather}
After rescaling the variable $\tau$ (i.e. after conjugating by the
element of $\PSL_2(\RR)$ represented by the matrix in
(\ref{eqn:series:conj}))
\begin{gather}\label{eqn:series:conj}
\left(
  \begin{array}{cc}
    1/2 & 0 \\
    0 & 1 \\
  \end{array}
\right)
\end{gather}
we obtain a space spanned by $(\eta(2\tau)/\eta(\tau))^{28}$ and its
reciprocal. This space is evidently preserved by the action of
$\Gamma_0(2)+2$. The Fricke involution $\tau\mapsto -1/2\tau$
interchanges (the rescalings of) $\tilde{F}^A_{2B}$ and
$\tilde{F}^{\forall}_{2B}$.

\subsubsection{$\gt{F}_{10B}$}

Elements of class $10B$ also lift to a unique class in the covering
group, and any such lift has trace $4$ on $\bigwedge(\gt{r})$.
Consulting Table~\ref{tab:series:Ru|frms} we see that
$\pm\hat{g}\sim 2^220^6/4^210^6$.
\begin{gather}
     \tilde{F}_{10B}^A(\tau)=\frac{4}{1}
          \frac{\eta(2\tau)^{2}\eta(20\tau)^6}
               {\eta(4\tau)^{2}\eta(10\tau)^6},\quad
     \tilde{F}_{10B}^{\forall}(\tau)=
          \frac{1}{4}
          \frac{\eta(4\tau)^{2}\eta(10\tau)^6}
          {\eta(2\tau)^{2}\eta(20\tau)^6}
\end{gather}
After rescaling we obtain the space spanned by
$\eta(\tau)^2\eta(10\tau)^6/\eta(2\tau)^2\eta(5\tau)^6$ and its
reciprocal. This space is preserved by the action of the group
$\Gamma_0(10)$.

\subsubsection{$\gt{F}_{14A}$, $\gt{F}_{14B}$, $\gt{F}_{14C}$}

Similar to the previous cases, elements of class $14A$, $14B$, or
$14C$, all lift to unique classes or order $28$ in the cover. Any
such lift has trace $4$ on $\bigwedge(\gt{r})$, and we have
$\pm\hat{g}\sim 28^4/14^4$.
\begin{gather}
     \tilde{F}_{14ABC}^A(\tau)=\frac{4}{1}
          \frac{\eta(28\tau)^{4}}{\eta(14\tau)^{4}},\quad
     \tilde{F}_{14ABC}^{\forall}(\tau)=\frac{1}{4}
          \frac{\eta(14\tau)^{4}}{\eta(28\tau)^{4}}
\end{gather}
Then similar to the $2B$ case we see that a conjugate of
$\Gamma_0(2)+2$ preserves the spaces $\gt{F}_{14ABC}$.

\subsubsection{$\gt{F}_{26A}$, $\gt{F}_{26B}$, $\gt{F}_{26C}$}

Again, all elements of class $26A$, $26B$, or $26C$, lift to unique
classes in the covering group, and have trace $4$ on
$\bigwedge(\gt{r})$. From Table~\ref{tab:series:Ru|frms} we have
$\pm\hat{g}\sim 4^252^2/2^226^2$, and this shows
\begin{gather}
     \tilde{F}_{26ABC}^A(\tau)=\frac{4}{1}
          \frac{\eta(4\tau)^{2}\eta(52\tau)^2}
               {\eta(2\tau)^{2}\eta(26\tau)^2},\quad
     \tilde{F}_{26ABC}^{\forall}(\tau)=
          \frac{1}{4}
          \frac{\eta(2\tau)^{2}\eta(26\tau)^2}
          {\eta(4\tau)^{2}\eta(52\tau)^2}.
\end{gather}
Here we obtain a representation of (the conjugate by
(\ref{eqn:series:conj}) of) the genus zero group $\Gamma_0(26)+26$.
The Fricke involution interchanges the two one dimensional
representations of the subgroup (conjugate to) $\Gamma_0(26)$.
(Incidentally, this subgroup is not a group of genus zero.)

\subsubsection{$\gt{F}_{29A}$, $\gt{F}_{29B}$}

Elements of class $29A$ or $29B$ have two lifts to the covering
group $2.\!\Ru$. A lift $\hat{g}$ of minimal order has order $29$,
and $-\hat{g}$ then has trace $29$ on $\bigwedge(\gt{r})$. The
corresponding Frame shape is $\hat{g}\sim 29^2/1^2$, and thus we
have
\begin{gather}
     \tilde{F}_{29AB}^A(\tau)=\frac{29}{1}
          \frac{\eta(29\tau)^2}{\eta(\tau)^2},\quad
     \tilde{F}_{29AB}^{\forall}(\tau)=
          \frac{1}{29}\frac{\eta(\tau)^2}{\eta(29\tau)^2}.
\end{gather}
The spaces $\gt{F}_{29AB}$ are each preserved by the action of
$\Gamma_0(29)+29$. The Fricke involution $\tau\mapsto -1/29\tau$
interchanges the one dimensional representations of $\Gamma_0(29)$
spanned by $\tilde{F}^A_{29A}$ and $\tilde{F}_{29A}^{\forall}$,
respectively (and similarly for the class $29B$).

\bigskip

We conclude with the following observation.
\begin{quote}
For each $g\in \Ru$ the functions $\tilde{F}^A_g(\tau)$ and
$\tilde{F}^{\forall}_g(\tau)$ span a (two dimensional)
representation of a discrete subgroup of $\PSL_2(\RR)$ that is
commensurable with $\PSL_2(\ZZ)$ and has genus zero.
\end{quote}

\comment{\begin{rmk} We remark here that the third group of Janko
(also a non-monstrous sporadic group) admits constructions by vertex
operators directly analogous to those given for the Rudvalis group
here and in the companion paper \cite{DunVARuI}. Consequently one
may define two dimensional spaces of functions $\gt{F}_g$ for each
$g$ in Janko's group $J_3$, and similar to the above, the spaces
$\gt{F}_g$ depend upon a choice of lift from $J_3$ to a covering
group $2\times 3.J_3$. We find once again that the lift may always
be chosen in such a way that $\gt{F}_g$ has the genus zero property
formulated above. These observations will receive a more thorough
treatment in the forthcoming article \cite{DunVAJ3}.
\end{rmk}}

\newcommand{\etalchar}[1]{$^{#1}$}

\end{document}